\documentclass{amsart}
\usepackage[alphabetic]{amsrefs}
\usepackage[all]{xy}
\usepackage[mathscr]{eucal}
\usepackage{graphicx,amssymb}
\CompileMatrices

\newtheorem{theorem}{Theorem}
\newtheorem{lemma}[theorem]{Lemma}
\newtheorem{proposition}[theorem]{Proposition}
\newtheorem{corollary}[theorem]{Corollary}

\theoremstyle{definition}
\newtheorem{definition}[theorem]{Definition}
\newtheorem{example}[theorem]{Example}
\newtheorem*{notation}{Notation}
\newtheorem*{observation}{Observation}

\theoremstyle{remark}
\newtheorem{remark}[theorem]{Remark}

\numberwithin{equation}{section}

\begin{document}

\newcommand{\Diff}{\operatorname{Diff}}
\newcommand{\Homeo}{\operatorname{Homeo}}
\newcommand{\Hom}{\operatorname{Hom}}
\newcommand{\Exp}{\operatorname{Exp}}
\newcommand{\Orb}{\operatorname{\textup{Orb}}}

\newcommand{\COrb}{\operatorname{\star\textup{Orb}}}
\newcommand{\CROrb}{\operatorname{\scriptscriptstyle{\blacklozenge}\scriptstyle\textup{Orb}}}
\newcommand{\ssslozenge}{\scriptscriptstyle{\blacklozenge}}
\newcommand{\ssstriangledown}{\scriptscriptstyle{\blacktriangledown}}

\newcommand{\Stwo}{\mbox{$\displaystyle S^2$}}
\newcommand{\Sn}{\mbox{$\displaystyle S^n$}}
\newcommand{\supp}{\operatorname{supp}}
\newcommand{\intr}{\operatorname{int}}
\newcommand{\kernel}{\operatorname{ker}} \newcommand{\A}{\mathbb{A}}
\newcommand{\B}{\mathbb{B}} \newcommand{\C}{\mathbb{C}}
\newcommand{\D}{\mathbb{D}} \newcommand{\E}{\mathbb{E}}
\newcommand{\F}{\mathbb{F}} \newcommand{\G}{\mathbb{G}}
\newcommand{\Hh}{\mathbb{H}} \newcommand{\I}{\mathbb{I}}
\newcommand{\J}{\mathbb{J}} \newcommand{\K}{\mathbb{K}}
\newcommand{\Ll}{\mathbb{L}} \newcommand{\M}{\mathbb{M}}
\newcommand{\N}{\mathbb{N}} \newcommand{\Oo}{\mathbb{O}}
\newcommand{\Pp}{\mathbb{P}} \newcommand{\Q}{\mathbb{Q}}
\newcommand{\R}{\mathbb{R}} \newcommand{\Ss}{\mathbb{S}}
\newcommand{\T}{\mathbb{T}} \newcommand{\U}{\mathbb{U}}
\newcommand{\V}{\mathbb{V}} \newcommand{\W}{\mathbb{W}}
\newcommand{\X}{\mathbb{X}} \newcommand{\Y}{\mathbb{Y}}
\newcommand{\Z}{\mathbb{Z}} \newcommand{\kk}{\mathbb{k}}
\newcommand{\orbify}[1]{\ensuremath{\mathcal{#1}}}
\newcommand{\starfunc}[1]{\ensuremath{{}_\star{#1}}}
\newcommand{\lozengefunc}[1]{\ensuremath{{}_{\scriptscriptstyle{\blacklozenge}}{#1}}}
\newcommand{\redfunc}[1]{\ensuremath{{}_\bullet{#1}}}

\newcommand{\OrbDiff}{\ensuremath{\Diff_{\Orb}}}
\newcommand{\RedOrbDiff}{\ensuremath{\Diff_{\textup{red}}}}
\newcommand{\COrbDiff}{\ensuremath{\Diff_{\COrb}}}
\newcommand{\CROrbDiff}{\ensuremath{\Diff_{\CROrb}}}
\newcommand{\OrbMaps}{\ensuremath{C_{\Orb}}}
\newcommand{\RedOrbMaps}{\ensuremath{C_{\textup{red}}}}
\newcommand{\COrbMaps}{\ensuremath{C_{\COrb}}}
\newcommand{\CROrbMaps}{\ensuremath{C_{\CROrb}}}
\newcommand{\Frechet}{Fr\'{e}chet\ }
\newcommand{\Frechetnospace}{Fr\'{e}chet}

\title[The Stratified Structure of Spaces of Smooth Orbifold Mappings]{The Stratified Structure of Spaces of Smooth Orbifold Mappings}

% Information for first author
\author{Joseph E. Borzellino}
% Address of record for the research reported here
\address{Department of Mathematics, California Polytechnic State
  University, 1 Grand Avenue, San Luis Obispo, California 93407}
% Current address \curraddr{Department of Mathematics and Statistics,
%   Case Western Reserve University, Cleveland, Ohio 43403}
\email{jborzell@calpoly.edu}
% \thanks will become a 1st page footnote.  \thanks{The first author
%   was supported in part by NSF Grant \#000000.}

% Information for second author
\author{Victor Brunsden} \address{Department of Mathematics and
  Statistics, Penn State Altoona, 3000 Ivyside Park, Altoona,
  Pennsylvania 16601} \email{vwb2@psu.edu}
% \thanks{Support information for the second author.}

% General info

\date{\today} \commby{Editor} 
\begin{abstract}
We consider four notions of maps between smooth $C^r$ orbifolds $\orbify{O}$, $\orbify{P}$ with $\orbify{O}$ compact (without boundary). We show that one of these notions is natural and necessary in order to uniquely define the notion of orbibundle pullback. For the notion of complete orbifold map, we show that the corresponding set of $C^r$ maps between $\orbify{O}$ and $\orbify{P}$ with the $C^r$ topology carries the structure of a smooth $C^\infty$ Banach ($r$ finite)/\Frechet ($r=\infty$) manifold. For the notion of complete reduced orbifold map, the corresponding set of $C^r$ maps between $\orbify{O}$ and $\orbify{P}$ with the $C^r$ topology carries the structure of a smooth $C^\infty$ Banach ($r$ finite)/\Frechet ($r=\infty$) orbifold. The remaining two notions carry a stratified structure: The $C^r$ orbifold maps between $\orbify{O}$ and $\orbify{P}$ is locally a stratified space with strata modeled on smooth $C^\infty$ Banach ($r$ finite)/\Frechet ($r=\infty$) manifolds while the set of $C^r$ reduced orbifold maps 
between $\orbify{O}$ and $\orbify{P}$ locally has the structure of a stratified space with strata modeled on smooth 
$C^\infty$ Banach 
($r$ finite)/\Frechet ($r=\infty$) orbifolds. Furthermore, we give the explicit relationship between these notions of orbifold map. Applying our results to the special case of orbifold diffeomorphism groups, we show they inherit the structure of $C^\infty$ Banach ($r$ finite)/\Frechet ($r=\infty$) manifolds. In fact, for $r$ finite they are topological groups, and for $r=\infty$ they are convenient \Frechet Lie groups.
\end{abstract}

\maketitle

\setcounter{tocdepth}{1}
\tableofcontents

\section{Introduction}\label{IntroSection}

A well-known result in the theory of differentiable dynamical systems states that the set of $C^r$ mappings
$C^r(M,N)$ between $C^\infty$ manifolds $M$ and $N$ with $M$ compact has the structure of a $C^\infty$ Banach manifold.
If $r=\infty$, $C^\infty(M,N)$ becomes a $C^\infty$ \Frechet manifold. The local model at $f\in C^r(M,N)$ is $\mathscr{D}^r(f^*(TN))$, the space of smooth sections of the pullback tangent bundle $f^*(TN)$ equipped with the $C^r$ topology. $\mathscr{D}^r(f^*(TN))$ is a separable Banach space for $1\le r < \infty$ and a separable \Frechet space for $r=\infty$. For reference, see
\cite{MR0203742}, \cite{MR0248880}, \cite{MR0649788}, \cite{MR546811},  or \cite{MR1445290}.

We wish to extend these results to the set of $C^r$ maps from a compact smooth orbifold $\orbify{O}$ (without boundary) to a smooth orbifold $\orbify{P}$. Interestingly, there are different useful notions of a $C^r$ map between orbifolds. In \cite{MR1926425}, we defined a notion of (unreduced) $C^r$ orbifold map and the notion of {\em reduced} orbifold map. In \cite{BB2006}, we clarified these notions and showed that for a compact orbifold $\orbify{O}$ (without boundary), both the group $\OrbDiff^r(\orbify{O})$ of orbifold diffeomorphisms and the group $\Diff^r_{\textup{red}}(\orbify{O})$ of reduced orbifold diffeomorphisms equipped with the $C^r$ topology, carry the topological structure of a Banach manifold for finite $r$ and \Frechet manifold for $r=\infty$. In fact, we showed that
$\Diff^r_{\textup{red}}(\orbify{O})$ is a finite quotient of $\OrbDiff^r(\orbify{O})$. While our notion of orbifold map is more general than the one that typically appears in the literature, for example \cite{MR2359514}, our notion of {\em reduced} orbifold map agrees with that book's definition~1.3 which is the definition that appears most often.

In order to extend the classical structure result for maps between manifolds to maps between orbifolds, and to generalize our results on the orbifold diffeomorphism group, we will introduce two additional notions of orbifold maps, the {\em complete orbifold maps} and the {\em complete reduced orbifold maps}. As simple examples will show, the notion of complete orbifold map is necessary to give a well-defined notion of pullback orbibundle. The need to be careful when defining pullback orbibundles was already noted in the work of Moerdijk and Pronk \cite{MR1466622} and Chen and Ruan \cite{MR1950941}.

To reconcile our definitions of orbifold map with the existing literature using Lie groupoid theoretic approach to orbifolds we make the following remarks. Our complete orbifold maps are essentially equivalent to the the groupoid homomorphisms of Moerdijk \cite{MR1950948}, however, they are independent of any particular groupoid structure one imposes on an orbifold and thus are more natural for the kinds of questions we address here.  Moreover, Chen's definition of orbifold map \cite{MR2263948} agrees with our notion of complete reduced orbifold map up to conjugation.

Lastly, we note that if $M$ and $N$ (as above) are, in addition, $\Gamma$-manifolds ($\Gamma$, a compact Lie group), then the space $C^r_\Gamma(M,N)$ of $C^r$ equivariant maps from $M$ to $N$ is a closed $C^\infty$ Banach submanifold of $C^r(M,N)$ \cite{MR0277850}. In \cite{BB2006}*{Example~3.10}, we observed that for a so-called good orbifold 
$\orbify{O}=M/\Gamma$, the orbifold diffeomorphism group $\OrbDiff^r(\orbify{O})$ is strictly larger than $\Diff^r_\Gamma(M)$, the $\Gamma$-equivariant diffeomorphism group of $M$. The relationship between the space of smooth orbifold maps between good orbifolds $\orbify{O}_i=M_i/\Gamma$, and the space of equivariant maps $C^r_\Gamma(M_1,M_2)$ will be the focus of a future investigation.

We assume the reader is familiar with the notion of smooth $C^r$ orbifolds, and although there are many nice references for this background material such as the recently published book \cite{MR2359514}, we will use our previous work \cite{BB2006} as our standard reference for notation and needed definitions. We should note, however, that our definition of orbifold is modeled on the definition in Thurston \cite{Thurston78} and that the orbifolds that concern us here are referred to as {\em classical effective orbifolds} in \cite{MR2359514}. More precisely, for our definition of orbifolds, isotropy actions are always effective and we allow for singularities of codimension one. For those notions for which the existing literature is not entirely consistent, we will provide explicit definitions. Our main result is the following

\begin{theorem}\label{MainTheorem}
  Let $r\ge 1$ and let $\orbify{O}, \orbify{P}$ be smooth $C^r$ orbifolds
  (without boundary) with $\orbify{O}$ compact. Denote by $\COrbMaps^r(\orbify{O},\orbify{P})$ the set
  of $C^r$ complete  orbifold maps between $\orbify{O}$ and $\orbify{P}$ equipped with the $C^r$ topology. Let $\starfunc f\in\COrbMaps^r(\orbify{O},\orbify{P})$.
  Then $\COrbMaps^r(\orbify{O},\orbify{P})$ is a smooth $C^\infty$ manifold modeled locally on the
  topological vector space $\mathscr{D}^r_{\Orb}(\starfunc f^*(T\orbify{P}))$ of $C^r$
  orbisections of the pullback tangent orbibundle of $\orbify{P}$ equipped with the $C^r$
  topology. This separable vector space is a Banach space if $1\le r <
  \infty$ and is a \Frechet space if $r = \infty$.
\end{theorem}

As corollaries of theorem~\ref{MainTheorem}, we are able to prove the following structure results for our different notions of orbifold map. For the complete reduced orbifold maps we have

\begin{corollary}\label{CompleteReducedOrbifoldMapStructure}
  Let $r\ge 1$ and let $\orbify{O}, \orbify{P}$ be as above. Denote by $\CROrbMaps^r(\orbify{O},\orbify{P})$ the set
  of complete reduced $C^r$ orbifold maps between $\orbify{O}$ and $\orbify{P}$ equipped with the $C^r$ topology inherited from
  $\COrbMaps^r(\orbify{O},\orbify{P})$ as a quotient space.
  Then $\CROrbMaps^r(\orbify{O},\orbify{P})$ carries the structure of a smooth $C^\infty$ 
  Banach ($r$ finite)/\Frechet ($r=\infty$) orbifold.
\end{corollary}

This result essentially recovers the result of Chen \cite{MR2263948} for $r$ finite, where the $C^r$ maps defined there are shown to have the structure of a smooth Banach orbifold. We have the following structure result for orbifold maps.

\begin{corollary}\label{OrbifoldMapStructure} Let $r\ge 1$ and let $\orbify{O}, \orbify{P}$ be as above. Denote by $\OrbMaps^r(\orbify{O},\orbify{P})$ the set of $C^r$ orbifold maps between $\orbify{O}$ and $\orbify{P}$ equipped with the $C^r$ topology (as defined in \cite{BB2006}).
Then $\OrbMaps^r(\orbify{O},\orbify{P})$ carries the topological structure of a stratified space with strata modeled on 
smooth $C^\infty$ Banach ($r$ finite)/\Frechet ($r=\infty$) manifolds.
\end{corollary}

In section~\ref{StratificationExampleSection}, we illustrate this phenomenon with a concrete example.
Finally, for the reduced orbifold maps, we conclude

\begin{corollary}\label{ReducedOrbifoldMapStructure} Let $r\ge 1$ and let $\orbify{O}, \orbify{P}$ be as above. Denote by $\RedOrbMaps^r(\orbify{O},\orbify{P})$ the set of $C^r$ reduced orbifold maps between $\orbify{O}$ and $\orbify{P}$ equipped with the $C^r$ topology as a quotient space.
Then $\RedOrbMaps^r(\orbify{O},\orbify{P})$ carries the topological structure of a stratified space with strata modeled on 
smooth $C^\infty$ Banach ($r$ finite)/\Frechet ($r=\infty$) orbifolds.
\end{corollary}

We would like to point out in each of the above results we are claiming, in part, the existence of a smooth structure modeled on Banach or \Frechet spaces. While much of the finite dimensional smooth manifold theory carries over to the Banach category, the lack of a general implicit function theorem in \Frechet spaces can cause significant difficulties \cite{MR656198}. In particular, there can be many inequivalent notions of differential calculus \cite{MR0440592}. For finite order differentiability, a strong argument can be made that the Lipschitz categories $\text{Lip}^r$ are better suited to questions of calculus than the more common $C^r$ category. For our purposes, however, we have chosen to use the $C^r$ category for finite order differentiability and for infinite order differentiability, we use the convenient calculus as detailed in the monographs \cites{MR961256, MR1471480}.

The paper is divided into the following sections: Section~\ref{OrbifoldMapSection} will define the four notions of orbifold map that we will be considering and how these notions are related. Section~\ref{FunctionSpaceTopologySection} defines the $C^r$ topology on $\COrbMaps^r(\orbify{O},\orbify{P})$ with $\orbify{O}$ compact and proves corollary~\ref{CompleteReducedOrbifoldMapStructure} assuming theorem~\ref{MainTheorem}. Section~\ref{DiffeosSection} applies our results to the special case of orbifold diffeomophisms. Section~\ref{StratificationExampleSection} provides explicit examples to show that non-orbifold structure stratifications naturally arise. Section~\ref{TangentOrbibundleSection} will construct the pullback orbibundle for a smooth complete orbifold map and illustrate the necessity to use complete orbifold maps in order to get a unique notion of pullback. Section~\ref{ProofOfMainTheoremSection} recalls some results about the exponential map on orbifolds and contains the proof of theorem~\ref{MainTheorem}. Section~\ref{StratifiedNeighborhoodsSection} is devoted to proofs of corollaries~\ref{OrbifoldMapStructure} and \ref{ReducedOrbifoldMapStructure}. In section~\ref{InfiniteDimAnalysis}, we collect the results of infinite-dimensional analysis that we need to substantiate our smoothness claims.

\section{Four Notions of Orbifold Map}\label{OrbifoldMapSection}
We now discuss four related definitions of maps between orbifolds. The first notion we will define is that of a complete orbifold map. It is distinguished from our previous notions of orbifold map and reduced orbifold map \cite{BB2006}*{Section 3} in that we are going to keep track of all defining data. In what follows we use the notation of \cite{BB2006}*{Section 2}.

\begin{definition}\label{CompleteOrbiMap}  
  A $C^0$ \emph{complete orbifold map} $(f,\{\tilde f_x\},\{\Theta_{f,x}\})$ between locally
  smooth orbifolds $\orbify{O}_1$ and $\orbify{O}_2$ consists of the
  following:
  \begin{enumerate}
  \item A continuous map $f:X_{\orbify{O}_1}\to X_{\orbify{O}_2}$ of
    the underlying topological spaces.
  \item For each $y\in {\mathcal{S}_x}$, a group homomorphism
    $\Theta_{f,y}:\Gamma_{\mathcal{S}_x}\to\Gamma_{f(y)}$.
  \item A $\Theta_{f,y}$-equivariant lift $\tilde f_y:\tilde
    U_y\subset\tilde{U}_{\mathcal{S}_x}\to\tilde V_{f(y)}$ where
    $(\tilde U_y,\Gamma_{\mathcal{S}_x}, \rho_y, \phi_y)$ is an
    orbifold chart at $y$ and $(\tilde V_{f(y)},\Gamma_{f(y)},
    \rho_{f(y)}, \phi_{f(y)})$ is an orbifold chart at $f(y)$.  That
    is, the following diagram commutes:
    \begin{equation*}
      \xymatrix{{\tilde U_y}\ar[rr]^{\tilde f_y}\ar[d]&&{\tilde V_{f(y)}}\ar[d]\\
        {\tilde U_y}/\Gamma_{\mathcal{S}_x}\ar[rr]^>>>>>>>>>>%
        {{\tilde f_y}/\Theta_{f,y}(\Gamma_{\mathcal{S}_x})}\ar[dd]&&{\tilde
          V_{f(y)}}/\Theta_{f,y}(\Gamma_{\mathcal{S}_x})\ar[d]\\
        &&{\tilde V_{f(y)}}/\Gamma_{f(y)}\ar[d]\\
        U_y\subset U_{\mathcal{S}_x}\ar[rr]^{f}&&V_{f(y)}
      }
    \end{equation*}
\begin{sloppypar}
  \item[($\star 4$)] (Equivalence) Two complete orbifold maps $(f,\{\tilde f_x\},\{\Theta_{f,x}\})$ and
    $(g,\{\tilde g_x\},\{\Theta_{g,x}\})$ are considered equivalent if for each
    $x\in\orbify{O}_1$, $\tilde f_x=\tilde g_x$ as germs and $\Theta_{f,x}=\Theta_{g,x}$. That is,
    there exists an orbifold chart $(\tilde U_x,\Gamma_x)$ at $x$ such
    that ${\tilde f_x}\vert_{\tilde U_x}={\tilde g_x}\vert_{\tilde
      U_x}$ and $\Theta_{f,x}=\Theta_{g,x}$. Note that this implies that $f=g$.
\end{sloppypar}
  \end{enumerate}
\end{definition}

\begin{definition}\label{CompleteOrbiMapSmooth} A complete orbifold map
  $f:\orbify{O}_1\to\orbify{O}_2$ of $C^r$ smooth orbifolds is {\em
    $C^r$ smooth} if each of the local lifts $\tilde f_x$ may be
  chosen to be $C^r$. Given two orbifolds $\orbify{O}_i$, $i = 1,2$, the set of $C^r$
complete orbifold maps from $\orbify{O}_1$ to $\orbify{O}_2$ will be denoted by
$C^r_{\COrb}(\orbify{O}_1, \orbify{O}_2)$. 
\end{definition}

If we replace ($\star 4$) in definition~\ref{CompleteOrbiMap} by

\begin{enumerate}
\begin{sloppypar}
\item[(4)] (Equivalence) Two complete orbifold maps $(f,\{\tilde f_x\},\{\Theta_{f,x}\})$ and
    $(g,\{\tilde g_x\},\{\Theta_{g,x}\})$ are considered equivalent if for each
    $x\in\orbify{O}_1$, $\tilde f_x=\tilde g_x$ as germs. That is,
    there exists an orbifold chart $(\tilde U_x,\Gamma_x)$ at $x$ such
    that ${\tilde f_x}\vert_{\tilde U_x}={\tilde g_x}\vert_{\tilde
      U_x}$ (which as before implies $f=g$),
\end{sloppypar}
\end{enumerate}
where we have dropped the requirement that $\Theta_{f,x}=\Theta_{g,x}$, we recover the notion of {\em orbifold map} $(f,\{\tilde f_x\})$ which appeared in \cite{BB2006}*{Section 3}. Thus, the set of orbifold maps $C^r_{\Orb}(\orbify{O}_1, \orbify{O}_2)$ can be regarded as the equivalence classes of complete orbifold maps under the less restrictive set-theoretic equivalence ($4$).
The following simple example is illustrative.

\begin{example}\label{RZ2Example} Let $\orbify{O}$ be the orbifold
  $\R/\Z_2$ where $\Z_2$ acts on $\R$ via $x \to -x$ and $f:\orbify{O}\to\orbify{O}$ is the
  constant map $f\equiv 0$. The underlying topological space
  $X_{\orbify{O}}$ of $\orbify{O}$ is $[0, \infty)$ and the isotropy
  subgoups are trivial for $x\in (0, \infty)$ and $\Z_2$ for $x = 0$.
  The map $\tilde f_0\equiv 0$ is a local equivariant lift of $f$ at
  $x=0$ using either of the homomorphisms $\Theta_{f,0}=\textup{Id}$
  or $\Theta'_{f,0}\equiv e$. Of course, for $x\ne 0$, we set $\tilde f_x\equiv 0$ and 
  $\Theta_{f,x}=\Theta'_{f,x}=\text{ the trivial homomorphism }\Gamma_x=e\mapsto e\in\Gamma_0=\Z_2$.
  Thus, as \emph{complete} orbifold maps $(f,\{\tilde f_x\},\{\Theta_{f,x}\})\ne (f,\{\tilde f_x\},\{\Theta'_{f,x}\})$.
  However, simply as orbifold maps, they are considered equal.
\end{example}

If we replace ($\star 4$) in definition~\ref{CompleteOrbiMap} by

\begin{enumerate}
\begin{sloppypar}
\item[($\scriptstyle{\blacklozenge}\displaystyle 4$)] (Equivalence) Two complete orbifold maps $(f,\{\tilde f_x\},\{\Theta_{f,x}\})$ and
    $(g,\{\tilde g_x\},\{\Theta_{g,x}\})$ are considered equivalent if $f=g$ and for each
    $x\in\orbify{O}_1$, we have $\Theta_{f,x}=\Theta_{g,x}$,
\end{sloppypar}
\end{enumerate}
where we have dropped the requirement that the germs of the lifts $\tilde f_x$ and $\tilde g_x$ agree, we obtain a new notion of orbifold map
$(f,\{\Theta_{f,x}\})$ which we call a \emph{complete reduced orbifold map}. The set of smooth complete reduced orbifold maps will be denoted by
$C^r_{\CROrb}(\orbify{O}_1, \orbify{O}_2)$. As before, it is clear that $C^r_{\CROrb}(\orbify{O}_1, \orbify{O}_2)$ is a set-theoretic quotient of $C^r_{\COrb}(\orbify{O}_1, \orbify{O}_2)$.

If we replace (4) in the definition of orbifold map, or ($\scriptstyle{\blacklozenge}\displaystyle 4$) in the definition of complete reduced orbifold map, by

\begin{enumerate}
\item[($\bullet 4$)] (Equivalence) Two orbifold maps $(f,\{\tilde f_x\})$ and
    $(g,\{\tilde g_x\})$, (or, complete reduced orbifold maps $(f,\{\Theta_{f,x}\})$ and
    $(g,\{\Theta_{g,x}\})$) are considered equivalent if $f=g$.
\end{enumerate}
we obtain the notion of {\em reduced orbifold map} from \cite{MR1926425}.
The set of smooth reduced orbifold maps will be denoted by
$\RedOrbMaps^r(\orbify{O}_1, \orbify{O}_2)$. Like before, it is clear that $\RedOrbMaps^r(\orbify{O}_1, \orbify{O}_2)$ is a set-theoretic quotient of both $C^r_{\Orb}(\orbify{O}_1, \orbify{O}_2)$ and $C^r_{\CROrb}(\orbify{O}_1, \orbify{O}_2)$.

\begin{notation} Since we will often need to distinguish between these various notions of orbifold maps, we will denote a complete orbifold map $(f,\{\tilde f_x\},\{\Theta_{f,x}\})$ by $\starfunc f$, and represent an orbifold map
$(f,\{\tilde f_x\})$ simply by $f$ as in \cite{BB2006}, a complete reduced orbifold map $(f,\{\Theta_{f,x}\})$ by $\lozengefunc f$, and a reduced orbifold map by $\redfunc f$.
\end{notation}

Diagrammatically, we have the following:
\begin{equation*}\label{MapDiagram}
      \xymatrix{&{\starfunc f\in C^r_{\COrb}(\orbify{O}_1, \orbify{O}_2)}%
      \ar[dl]_{q_{\scriptscriptstyle\blacklozenge}}\ar[dr]^{q}\ar@{-->}[dd]^{q_\star}&\\
        {\lozengefunc f\in C^r_{\CROrb}(\orbify{O}_1, \orbify{O}_2)}\ar[rd]^{q_{\scriptscriptstyle\blacktriangledown}} && {f\in C^r_{\Orb}(\orbify{O}_1, \orbify{O}_2)}\ar[ld]_{q_\bullet}\\
        &{\redfunc f\in \RedOrbMaps^r(\orbify{O}_1, \orbify{O}_2)} &
      }
\end{equation*}
where the $q$'s represent the respective set-theoretic quotient maps.
Understanding how these notions are related in the special case of the identity map is crucial in what follows.

\begin{example}(Lifts of the Identity Map)\label{IdentityMap} Consider the identity map $\textup{Id}:\orbify{O}\to\orbify{O}$. Let 
  $x\in\orbify{O}$ and $(\tilde U_x,\Gamma_x)$ be an orbifold chart at
  $x$. From the definition of orbifold map, it follows (since
  $\Gamma_x$ is finite) that there exists $\gamma\in\Gamma_x$ such
  that a lift $\widetilde{\textup{Id}}_x:\tilde U_x\to\tilde U_x$ is
  given by $\widetilde{\textup{Id}}_x(\tilde y)=\gamma\cdot\tilde y$
  for all $\tilde y\in\tilde U_x$. Since $\widetilde{\textup{Id}}_x$
  is $\Theta_{\textup{Id},x}$ equivariant we have for
  $\delta\in\Gamma_x$:
  \begin{alignat*}{2}\widetilde{\textup{Id}}_x(\delta\cdot\tilde y) &
    = %
\Theta_{\textup{Id},x}(\delta)\cdot\widetilde{\textup{Id}}_x(\tilde y) &\quad & \text{hence }\\
\gamma\delta\cdot\tilde y & = \Theta_{\textup{Id},x}(\delta)\gamma\cdot\tilde y & & \text{which implies}\\
& & & \text{since $\Gamma_x$ acts effectively that}\\
\gamma\delta & = \Theta_{\textup{Id},x}(\delta)\gamma & & \text{or, equivalently,}\\
\Theta_{\textup{Id},x}(\delta) & = \gamma\delta\gamma^{-1}
\end{alignat*}
Thus, the isomorphism $\Theta_{\textup{Id},x}$ is completely determined by the choice of local lift $\widetilde{\textup{Id}}_x$.
This implies that the group $\mathscr{ID}$ of orbifold maps covering the identity may be regarded as the same as the group $\starfunc{\mathscr{ID}}$ of complete orbifold maps covering the identity. That is, we have the bijective correspondence
$$(\textup{Id},\{\tilde y\mapsto\gamma_{x}\cdot\tilde y\},\{\Theta_{\textup{Id},x}\})\longleftrightarrow%
(\textup{Id},\{\tilde y\mapsto\gamma_{x}\cdot\tilde y\}).$$
Suppose now that $\{U_{x_i}\}$ is a countable (possibly finite) cover of $\orbify{O}$ by charts. Then $\mathscr{ID}$ can be regarded as a subgroup of the product $\prod \Gamma_{x_i}$ as in the proof of corollary~1.2 in \cite{BB2006}.
Two inner automorphisms,
$\delta\mapsto\gamma_i\delta\gamma_i^{-1}$, give rise to the same
automorphism of $\Gamma_x$ precisely when $\gamma_1=\zeta\gamma_2$
where $\zeta\in C(\Gamma_x)$, the center of $\Gamma_x$. Thus, if we let 
$C=C(\mathscr{ID})\subset\prod C(\Gamma_{x_i})$, then one can see that
the complete reduced lifts of the identity $\lozengefunc{\mathscr{ID}}\cong\starfunc{\mathscr{ID}}/C$, where the free $C$-action on
$\starfunc{\mathscr{ID}}$ is defined by 
$$(\zeta_i)\cdot (\textup{Id},\{\tilde y\mapsto\gamma_{x_i}\cdot\tilde y\},\{\Theta_{\textup{Id},x_i}\})=%
(\textup{Id},\{\tilde y\mapsto (\zeta_i\gamma_{x_i})\cdot\tilde y\},\{\zeta_i\Theta_{\textup{Id},x_i}\zeta^{-1}_i=\Theta_{\textup{Id},x_i}\}).$$
Also, note that the correspondence $\starfunc{\mathscr{ID}}\leftrightarrow\mathscr{ID}$ gives an isomorphism
$\lozengefunc{\mathscr{ID}}\cong\mathscr{ID}/C$ which in turn is isomorphic to $\textup{Inn}(\mathscr{ID})$, the group of inner automorphisms of $\mathscr{ID}$. Thus, we have the exact sequence
$$1\rightarrow C(\mathscr{ID})\rightarrow\starfunc{\mathscr{ID}}=\mathscr{ID}\rightarrow\lozengefunc{\mathscr{ID}}\rightarrow 1$$
\end{example}

\begin{notation} For a (not necessarily compact) orbifold $\orbify{N}$, we will use the notation $\mathscr{ID}_{\orbify{N}}$ to denote the group of orbifold lifts of the identity map $\text{Id}:\orbify{N}\to\orbify{N}$. Suppose $f:\orbify{O}_1\to\orbify{O}_2$ and let the orbifold $\orbify{N}$ be an open neighborhood of the image
$f(\orbify{O}_1)$. For an orbifold map ${}_{\{\cdot\}} f$ (of any type) and $I=(\text{Id},\{\eta_x\cdot\tilde y\})\in\mathscr{ID}_{\orbify{N}}$ 
we can compute $I\circ {}_{\{\cdot\}} f$. Namely,
\begin{align*}
I\circ\starfunc f&  =  I\circ (f,\{\tilde f_x\},\{\Theta_{f,x}\})= (f, \{\eta_x\cdot\tilde f_x\}, %
\{\gamma\mapsto\eta_x\Theta_{f,x}(\gamma)\eta_x^{-1}\})\\
I\circ\lozengefunc f&  =  I\circ (f,\{\Theta_{f,x}\})= (f, \{\gamma\mapsto\eta_x\Theta_{f,x}(\gamma)\eta_x^{-1}\})\\
I\circ f&  =  I\circ (f,\{\tilde f_x\})= (f, \{\eta_x\cdot\tilde f_x\})\\
I\circ\redfunc f& =  \redfunc f
\end{align*}
Suppose $\{\Gamma_x\}$ denotes the family of isotropy groups for an orbifold 
$\orbify{N}$ and for subgroups $\Lambda_x\subset\Gamma_x$, let $\{\Lambda_x\}$ denote the corresponding family of subgroups. In what follows, we will use the notation $(\mathscr{ID}_{\orbify{N}})_{\{\Lambda_x\}}$ for the subgroup of $\mathscr{ID}_{\orbify{N}}$ defined by
$$\{I\in\mathscr{ID}_{\orbify{N}}\mid I=(\text{Id},\{\tilde y\mapsto\lambda_x\cdot\tilde y\})%
\text{ where }\lambda_x\in\Lambda_x\text{ for all }x\}.$$
Lastly, for a fixed orbifold map ${}_{\{\cdot\}} f$ (of any type), we let 
$(\mathscr{ID}_{\orbify{N}})\cdot {}_{\{\cdot\}} f$ denote the orbit under the action of $\mathscr{ID}_{\orbify{N}}$:
$$(\mathscr{ID}_{\orbify{N}})\cdot {}_{\{\cdot\}} f=\{I\circ {}_{\{\cdot\}} f\mid I\in\mathscr{ID}_{\orbify{N}}\}$$
and
we let $(\mathscr{ID}_{\orbify{N}})_{{}_{\{\cdot\}} f}$ denote the corresponding isotropy subgroup of ${}_{\{\cdot\}} f$ under the action of $\mathscr{ID}_{\orbify{N}}$:
$$\{I\in\mathscr{ID}_{\orbify{N}}\mid I\circ {}_{\{\cdot\}} f={}_{\{\cdot\}} f\}.$$

It is also important to note that $\mathscr{ID}_{\orbify{N}}$ is a finite group in the special case that the source orbifold $\orbify{O}_1$ is compact: one may choose the open neighborhood 
$\orbify{N}$ of $f(\orbify{O}_1)$ to be relatively compact and since $\orbify{N}$ can be covered by finitely many orbifold charts $\{U_{x_i}\}$, the observation that $\mathscr{ID}\subset\prod\Gamma_{x_i}$ from example~\ref{IdentityMap} is enough to show that, in this case, $\mathscr{ID}_{\orbify{N}}$ is finite.

\end{notation}

\subsection*{Implications for the definition of orbifold structure}
Recall the following commutative diagram of maps which appears in the definition of a smooth classical effective orbifold \cite{BB2006}:

\begin{equation*}
    \xymatrix{{\tilde U_z}\ar[rr]^{\tilde\psi_{zx}}\ar[d]&&{\tilde U_x}\ar[d]\\
      {U_z\cong\tilde U_z/\Gamma_z}\ar[rr]^{\psi_{zx}}&&{\tilde U_x/\Gamma_x\cong U_x}\\
     }
  \end{equation*}
where for
  a neighborhood $U_z\subset U_x$ with corresponding $\tilde U_z$, and isotropy group $\Gamma_z$, there is an open embedding $\tilde\psi_{zx}:\tilde
  U_z \to \tilde U_x$ covering the inclusion $\psi_{zx}:U_z\hookrightarrow U_x$ and an injective homomorphism
  $\theta_{zx}:\Gamma_z \to \Gamma_x$ so that $\tilde\psi_{zx}$ is
  equivariant with respect to $\theta_{zx}$. For the standard definition of orbifold which appears in the literature, it is understood that
  $\tilde{\psi}_{zx}$ is defined only up to composition with elements of $\Gamma_x$, and $\theta_{zx}$ defined only up to conjugation by elements of 
  $\Gamma_x$. However, here, we may regard $\psi_{zx}$ as being from any of the notions of orbifold map we have defined, thus giving an orbifold 
  $\orbify{O}$, or more precisely, an orbifold atlas for $\orbify{O}$, one of four different structures depending on how one keeps track of lifts 
  $\psi_{zx}$ and homomorphisms $\theta_{zx}$. Thus, it makes sense to speak of a {\em complete orbifold structure} $\starfunc{\orbify{O}}$,  a 
  {\em complete reduced orbifold structure} $\lozengefunc{\orbify{O}}$, an {\em orbifold structure} $\orbify{O}$, and lastly,  
  a {\em reduced orbifold structure} $\redfunc{\orbify{O}}$. Thus, the standard definition of orbifold would correspond to our notion of a reduced orbifold structure. The reader should take care to note that the term {\em reduced orbifold} also has been used in the study of so-called noneffective orbifolds \cite{MR1950941}. Our use of the term {\em reduced orbifold structure} is unrelated to this. In this paper, the term orbifold will require that the chart maps $\psi_{zx}$ be regarded as orbifold maps in $C_{\Orb}^r(U_z,U_x)$ as defined above.  We also point out that there is no fundamental difference between a complete orbifold structure $\starfunc{\orbify{O}}$ and an orbifold structure $\orbify{O}$ and that any reduced orbifold structure $\redfunc{\orbify{O}}$ is obtained as a quotient an orbifold structure $\orbify{O}$ by the action of $\mathscr{ID}$ on orbifold atlases. This follows from example~\ref{IdentityMap} and the fact that any two lifts of $\psi_{zx}$ must differ by a lift of the identity map on $U_x$. Lastly, we remark that, in general, for an orbifold structure $\orbify{O}$, $\tilde{\psi}_{zx}\ne\tilde{\psi}_{yx}\circ\tilde{\psi}_{zy}$ when
  $U_z\subset U_y\subset U_x$, but there will be an element
  $\delta\in\Gamma_x$ such that
  $\delta\cdot\tilde{\psi}_{zx}=\tilde{\psi}_{yx}\circ\tilde{\psi}_{zy}$
  and $\delta\cdot
  \theta_{zx}(\gamma)\cdot\delta^{-1}=\theta_{yx}\circ\theta_{zy}(\gamma)$.

\subsection*{Relationship among the different notions of orbifold map}
In this subsection we give a series of lemmas that discuss the relationship among the various notions of orbifold map for a fixed map $f:\orbify{O}_1\to\orbify{O}_2$. In section~\ref{FunctionSpaceTopologySection}, we will topologize these sets of mappings and discuss the local structure of these relationships. Our first lemma makes explicit the relationship between the complete reduced orbifold maps and the complete orbifold maps.

\begin{lemma}\label{PointwiseCOrbQuotientCROrb}
Let $\starfunc f,\starfunc f'\in C^r_{\COrb}(\orbify{O}_1, \orbify{O}_2)$ be complete orbifold maps which represent the same complete reduced orbifold map. That is, $\lozengefunc f=\lozengefunc f'\in C^r_{\CROrb}(\orbify{O}_1, \orbify{O}_2)$, so that 
$\starfunc f=(f,\{\tilde f_x\},\{\Theta_{f,x}\})$ and $\starfunc f'=(f,\{\tilde f'_x\},\{\Theta_{f,x}\})$. Let 
$C_x=C_{\Gamma_{f(x)}}(\Theta_{f,x}(\Gamma_x))$ denote the centralizer of $\Theta_{f,x}(\Gamma_x)$ in $\Gamma_{f(x)}$. Then there is an orbifold 
$\orbify{N}$ which is an open neighborhood of $f(\orbify{O}_1)$ in $\orbify{O}_2$ and an orbifold map 
$I\in (\mathscr{ID}_{\orbify{N}})_{\{C_x\}}$ such that $\starfunc f=I\circ\starfunc f'$. Moreover, if the stated condition holds for two complete orbifold maps $\starfunc f$ and $\starfunc f'$, then $\lozengefunc f=\lozengefunc f'$.
\end{lemma}

\begin{proof} Since $\tilde f_x$ and $\tilde f'_x$ are local lifts of the same map $f$, there exists $\eta_x\in\Gamma_{f(x)}$ such that
$\tilde f'_x(\tilde y)=\eta_x\cdot \tilde f_x(\tilde y)$ for all $\tilde y\in \tilde U_x$. Thus, for all $\gamma\in\Gamma_x$ we have,
on one hand, the equivariance relation $\tilde f'_x(\gamma\cdot\tilde y)=\Theta_{f,x}(\gamma)\cdot\tilde f'_x(\tilde y)$ while on the other hand, the equivariance relation must be
$\tilde f'_x(\gamma\cdot\tilde y)=\eta_x\cdot\tilde f_x(\gamma\cdot\tilde y)=\eta_x\Theta_{f,x}(\gamma)\cdot\tilde f_x(\tilde y)=
\eta_x\Theta_{f,x}(\gamma)\eta_x^{-1}\cdot\tilde f'_x(\tilde y)$. This implies that $\Theta_{f,x}(\gamma)=\eta_x\Theta_{f,x}(\gamma)\eta_x^{-1}$ and thus $\eta_x\in C_x$. The orbifold $\orbify{N}$ may be taken to be 
$\orbify{N}=\cup_{x\in\orbify{O}_1}V_{f(x)}$, where $V_{f(x)}$ is an orbifold chart about $f(x)\in\orbify{O}_2$. We have thus shown the first statement of the lemma, and the last statement is clear from our computation above and the definitions.
\end{proof}

\begin{example}\label{PointwiseCOrbQuotientCROrbNotFree} Let $\orbify{O}$ be as in example~\ref{RZ2Example}. Consider the complete orbifold map 
$\starfunc f=(f,\{\tilde f_x\},\{\Theta_{f,x}\})$ 
which covers the inclusion map 
$f:\orbify{O}\to\orbify{O}\times\orbify{O}\times\orbify{O}$, $y\mapsto (y,0,0)$, where $\tilde f_x(\tilde y)=(\tilde y,0,0)$ and $\Theta_{f,0}(\gamma)=(\gamma,e,e)\in\Gamma_{(0,0,0)}=\Z_2\times\Z_2\times\Z_2$ (for $x\ne 0$, $\Theta_{f,x}$ is the trivial homomorphism since $\Gamma_x=\{e\}$). Now, 
$\eta_0=(e,\gamma,e)\in C_{\Gamma_{(0,0,0)}}(\Theta_{f,0}(\Gamma_0))$ and 
$\eta_0\cdot\tilde f_0(\tilde y)=\eta_0\cdot (\tilde y,0,0)=(\tilde y,0,0)=\tilde f_0(\tilde y)$. Now let the (finite) group $(\mathscr{ID}_\orbify{N})_{\{C_x\}}$ be as in lemma~\ref{PointwiseCOrbQuotientCROrb}. For fixed $\starfunc f$, 
let $(\mathscr{ID}_\orbify{N})_{\{C_x\}}\cdot\starfunc f$ denote the orbit of $\starfunc f$. This example shows that the orbit map 
$(\textup{Id},\{\lambda_x\cdot\tilde z\})\mapsto (\textup{Id},\{\lambda_x\cdot\tilde z\})\circ\starfunc f$  may have nontrivial, (but finite) isotropy. Thus, if we let 
$(\mathscr{ID}_\orbify{N})_{\starfunc f}\subset(\mathscr{ID}_\orbify{N})_{\{C_x\}}$ denote the isotropy subgroup of $\starfunc f$, then 
$(\mathscr{ID}_\orbify{N})_{\{C_x\}}/(\mathscr{ID}_\orbify{N})_{\starfunc f}\cong (\mathscr{ID}_\orbify{N})_{\{C_x\}}\cdot\starfunc f$ is a homeomorphism (of discrete sets).
\end{example}

The next lemma describes the relationship between the orbifold maps and the reduced orbifold maps.

\begin{lemma}\label{PointwiseOrbQuotientRedOrb}
Let $f,f'\in \OrbMaps^r(\orbify{O}_1, \orbify{O}_2)$ be orbifold maps which represent the same reduced orbifold map. That is, $\redfunc f=\redfunc f'\in \RedOrbMaps^r(\orbify{O}_1, \orbify{O}_2)$, so that 
$f=(f,\{\tilde f_x\})$ and $f'=(f,\{\tilde f'_x\})$. Then there is an orbifold 
$\orbify{N}$ which is an open neighborhood of $f(\orbify{O}_1)$ in $\orbify{O}_2$ and an orbifold map 
$I=(\textup{Id},\{\eta_x\cdot\tilde z\})\in\mathscr{ID}_{\orbify{N}}$ with 
$\eta_x\in\Gamma_{f(x)}$ such that $f=I\circ f'$. Moreover, if the stated condition holds for two
orbifold maps $f$ and $f'$, then $\redfunc f=\redfunc f'$.
\end{lemma}

\begin{proof} $\orbify{N}$ can be chosen as in lemma~\ref{PointwiseCOrbQuotientCROrb}, and the proof follows from corollary~1.2 in \cite{BB2006}.
\end{proof}

\begin{remark}\label{PointwiseOrbQuotientRedOrbIsNotFree} Similar to the situation described in example~\ref{PointwiseCOrbQuotientCROrbNotFree}, example~\ref{RZ2Example} shows that the orbit map $\mathscr{ID}_{\orbify{N}}\cdot f$ may have nontrivial isotropy.  
\end{remark}

Next, we describe the relationship between the complete orbifold maps and the orbifold maps.

\begin{lemma}\label{PointwiseCOrbQuotientOrb}
Let $\starfunc f,\starfunc f'\in \COrbMaps^r(\orbify{O}_1, \orbify{O}_2)$ be complete orbifold maps which represent the same orbifold map. That is, $f=f'\in \OrbMaps^r(\orbify{O}_1, \orbify{O}_2)$, so that 
$\starfunc f=(f,\{\tilde f_x\},\{\Theta_{f,x}\})$ and $\starfunc f'=(f,\{\tilde f_x\},\{\Theta'_{f,x}\})$. Then, 
for each $x\in\orbify{O}_1$, $\gamma\in\Gamma_x$, and $\tilde y\in\tilde U_x$ we have
$$\left\{{\left[\Theta'_{f,x}(\gamma)\right]}^{-1}\left[\Theta_{f,x}(\gamma)\right]\right\}\cdot %
\tilde f_x(\tilde y)=\tilde f_x(\tilde y).$$
Moreover, if the stated condition holds for two
complete orbifold maps $\starfunc f$ and $\starfunc f'$, then $f=f'$.
\end{lemma}

\begin{proof} For all $\gamma\in\Gamma_x$ and $\tilde y\in\tilde U_x$, we have, $\Theta_{f,x}(\gamma)\cdot\tilde f_x(\tilde y)=\tilde f_x(\gamma\cdot\tilde y)=\Theta'_{f,x}(\gamma)\cdot\tilde f_x(\tilde y)$ and the first statement follows. To see the last statement, let
$\starfunc f=(f,\{\tilde f_x\},\{\Theta_{f,x}\})$ and $\starfunc f'=(f,\{\tilde f'_x\},\{\Theta'_{f,x}\})$. Note that the condition stated implies that 
$\tilde f_x(\tilde y)=\Theta_{f,x}(e)\cdot\tilde f_x(\tilde y)=\Theta'_{f,x}(e)\cdot\tilde f'_x(\tilde y)=\tilde f'_x(\tilde y)$.
\end{proof}

\begin{remark} Notice that this relationship is qualitatively different than the relationships described in
lemmas~\ref{PointwiseCOrbQuotientCROrb},
\ref{PointwiseOrbQuotientRedOrb} and \ref{PointwiseCROrbQuotientRedOrb}, in that it is given as an equality of \emph{actions}
of $\Theta_{f,x}(\Gamma_x)$, $\Theta'_{f,x}(\Gamma_x)$ on the image $\tilde f_x(\tilde U_x)$ and not as an equality of $\Theta_{f,x}$ and 
$\Theta'_{f,x}$ as homomorphisms  themselves. That is, the representation of $\Theta_{f,x}(\Gamma_x)$ and $\Theta'_{f,x}(\Gamma_x)$ induce actions that when restricted to $\tilde f_x(\tilde U_x)$ are equal.
\end{remark}

\begin{remark} Example~\ref{RZ2Example} exhibits the behavior described in lemma~\ref{PointwiseCOrbQuotientOrb}. A slightly less trivial example is to consider the inclusion map of example~\ref{PointwiseCOrbQuotientCROrbNotFree}: $f:\orbify{O}\to\orbify{O}\times\orbify{O}\times\orbify{O}$, 
$y\mapsto (y,0,0)$, where $\tilde f_x(\tilde y)=(\tilde y,0,0)$. Note that $\tilde f_0$ is equivariant with respect to both 
$\Theta_{f,0}(\gamma)=(\gamma,e,e)$ and $\Theta'_{f,0}(\gamma)=(\gamma,\gamma,\gamma)$.
\end{remark}

The next two lemmas describe the relationship between the complete reduced orbifold maps and the reduced orbifold maps. 
Given the conclusion of corollary~\ref{ReducedOrbifoldMapStructure}, this relationship is necessarily more complicated.

\begin{lemma}\label{PointwiseCROrbQuotientRedOrb}
Let $\lozengefunc f\in \CROrbMaps^r(\orbify{O}_1, \orbify{O}_2)$ be a complete reduced orbifold map and let 
$\orbify{N}$ be an open neighborhood of $f(\orbify{O}_1)$ in $\orbify{O}_2$. Let 
$I=(\textup{Id},\{\eta_x\cdot\tilde z\})\in\mathscr{ID}_{\orbify{N}}$ with 
$\eta_x\in\Gamma_{f(x)}$. If the complete reduced orbifold map $\lozengefunc f'=I\circ\lozengefunc f$, then 
$\redfunc f=\redfunc f'\in \RedOrbMaps^r(\orbify{O}_1, \orbify{O}_2)$. Furthermore, if $\lozengefunc f'$ is in the orbit
$(\mathscr{ID}_{\orbify{N}})\cdot\lozengefunc f$, then $\Theta'_{f,x}(\gamma)=\eta_x\Theta_{f,x}(\gamma)\eta_x^{-1}$.
\end{lemma}

\begin{proof} Let $\lozengefunc f=(f,\{\Theta_{f,x}\})$ and $\lozengefunc f'=(f',\{\Theta'_{f,x}\})$. Let $\tilde f_x,\tilde f'_x$ be local lifts equivariant with respect to $\Theta_{f,x},\Theta'_{f,x}$ respectively. 
Since $\lozengefunc f'=I\circ\lozengefunc f$, $\tilde f'_x(\tilde y)=\eta_x\cdot \tilde f_x(\tilde y)$ for all $\tilde y\in \tilde U_x$, so
$\tilde f_x$ and $\tilde f'_x$ are local lifts of the same map $f=f'$. This implies $\redfunc f=\redfunc f'$.
The last statement follows from the way $\mathscr{ID}_{\orbify{N}}$ acts on $\lozengefunc f$.
\end{proof}

\begin{remark} Here, like before, the orbit map 
$(\textup{Id},\{\eta_x\cdot\tilde z\})\mapsto (\textup{Id},\{\eta_x\cdot\tilde z\})\circ\lozengefunc f$  may have nontrivial, (but finite) isotropy. In fact,
$(\mathscr{ID}_\orbify{N})_{\lozengefunc f}=(\mathscr{ID}_\orbify{N})_{\{C_x\}}$, the orbifold map lifts of the identity given by elements of $C_{\Gamma_{f(x)}}(\Theta_{f,x}(\Gamma_x))$ described in lemma~\ref{PointwiseCOrbQuotientCROrb}.
\end{remark}

In light of lemma~\ref{PointwiseCROrbQuotientRedOrb}, we define an equivalence relation the preimage 
$q^{-1}_{\ssstriangledown}(\redfunc f)$:
\begin{equation}
\tag{$\dagger$}\lozengefunc f=(f,\{\Theta_{f,x}\})\sim\lozengefunc f'=(f',\{\Theta'_{f,x}\})\Longleftrightarrow \lozengefunc f'=I\circ\lozengefunc f
\end{equation}
for some $I\in\mathscr{ID}_\orbify{N}$. That is, for all $x\in\orbify{O}_1$, there exists $\eta_x\in\Gamma_{f(x)}$ such that 
$\Theta'_{f,x}(\gamma)=\eta_x\Theta_{f,x}(\gamma)\eta_x^{-1}$ for all $\gamma\in\Gamma_x$. Denote the equivalence class of
$\lozengefunc f$ by $[\lozengefunc f]$.

\begin{lemma}\label{PointwiseConjCROrbQuotientRedOrb}
Let $[\lozengefunc f=(f,\{\Theta_{f,x}\})]\ne [\lozengefunc f'=(f,\{\Theta'_{f,x}\})]$ be different equivalence classes of complete reduced orbifold maps which represent the same reduced orbifold map. That is, $\redfunc f=\redfunc f'\in \RedOrbMaps^r(\orbify{O}_1, \orbify{O}_2)$. Then there exist local lifts $\{\tilde f_x\}$ which are equivariant with respect to both $\{\Theta_{f,x}\}$ and 
$\{\eta_x\Theta'_{f,x}\eta_x^{-1}\}$ with $\{\Theta_{f,x}\}\ne \{\eta_x\Theta'_{f,x}\eta_x^{-1}\}$ as homomorphisms.
However, for each $x\in\orbify{O}_1$, $\gamma\in\Gamma_x$, and $\tilde y\in\tilde U_x$ we have as actions on $\tilde f_x(\tilde U_x)$
\begin{equation}
\Theta_{f,x}(\gamma)\cdot\tilde f_x(\tilde y)=\eta_x\Theta'_{f,x}(\gamma)\eta_x^{-1}\cdot\tilde f_x(\tilde y).
\tag{$\ddagger$}
\end{equation}
\end{lemma}
\begin{proof} Since $\redfunc f=\redfunc f'$, there exists $\eta_x\in\Gamma_{f(x)}$ such that 
$\tilde f_x(\tilde y)=\eta_x\cdot\tilde f'_x(\tilde y)$ for all $\tilde y\in\tilde U_x$. Thus, we conclude that 
$\tilde f_x=\eta_x\cdot\tilde f'_x$ is also equivariant with respect to $\eta_x\Theta'_{f,x}\eta_x^{-1}$. Since 
$[\lozengefunc f]\ne [\lozengefunc f']$, we have $\{\Theta_{f,x}\}\ne \{\eta_x\Theta'_{f,x}\eta_x^{-1}\}$ as homomorphisms.
\end{proof}

\begin{remark}\label{qtriangleFactors} Example~\ref{RZ2Example} illustrates the phenomena dealt with in lemma~\ref{PointwiseConjCROrbQuotientRedOrb}. 
Lemmas~\ref{PointwiseCROrbQuotientRedOrb} and \ref{PointwiseConjCROrbQuotientRedOrb} show that the quotient map 
$q_{\ssstriangledown}:\CROrbMaps^r(\orbify{O}_1,\orbify{O}_2)\to\RedOrbMaps^r(\orbify{O}_1,\orbify{O}_2)$ factors
$q_{\ssstriangledown}=q_{\ddagger}\circ q_{\dagger}$:
$$\xymatrix@1{
\lozengefunc f\ar@/_1.5pc/[rrrr]_{q_{\ssstriangledown}}\ar[rr]^{q_{\dagger}}   & &  [\lozengefunc f]\ar[rr]^{q_{\ddagger}} &  & \redfunc f
}$$
where $q_{\dagger}$, $q_{\ddagger}$ represent the quotient maps under the equivalences $(\dagger)$ and $(\ddagger)$, respectively.
\end{remark}

\section{Function Space Topologies}\label{FunctionSpaceTopologySection} It is easy to define a $C^s$ topology $(1\le s\le r)$ on the set of smooth complete orbifold maps $C^r_{\COrb}(\orbify{O},\orbify{P})$ with $\orbify{O}$ compact. Although much of what we do applies to noncompact $\orbify{O}$ we will assume $\orbify{O}$ to be compact. As such, implicit in some of the discussion is that $\orbify{O}$ has been equipped with a finite covering by orbifold charts. The topologies we define have already been shown to be independent of these choices of charts \cite{BB2006}.

\begin{definition}\label{CrTopOnComplete}Let $\starfunc f=(f,\{\tilde f_x\},\{\Theta_{f,x}\})$, 
$\starfunc g=(g,\{\tilde g_x\},\{\Theta_{g,x}\})\in\COrbMaps^r(\orbify{O},\orbify{P})$. Then a $C^s$ neighborhood of $\starfunc f$ is defined to be
\begin{multline*}
\mathscr{N}^s(\starfunc f,\varepsilon)=\{\starfunc g=(g,\{\tilde g_x\},\{\Theta_{g,x}\})\mid\\%
 g\in\mathscr{N}^s(f,\varepsilon)\text{ and }%
\theta_{f(x)z}\circ\Theta_{f,x}=\theta_{g(x)z}\circ\Theta_{g,x}\text{ for all }x\in\orbify{O}\}
\end{multline*}
where $\mathscr{N}^s(f,\varepsilon)$ is the $C^s$ orbifold map neighborhood of $f$ defined in \cite{BB2006}. 
$\theta_{f(x)z}\circ\Theta_{f,x}=\theta_{g(x)z}\circ\Theta_{g,x}$ is to be interpreted as follows: There is a small enough orbifold chart 
$\tilde U_x$ about $x$, such that the images of both $\tilde f_x(\tilde U_x)$ and
$\tilde g_x(\tilde U_x)$ are contained in a single orbifold chart $\tilde V_{z}$ and
$\theta_{f(x)z}\circ\Theta_{f,x}=\theta_{g(x)z}\circ\Theta_{g,x}$ where 
$\theta_{f(x)z},\theta_{g(x)z}:\Gamma_{f(x)},\Gamma_{g(x)}\hookrightarrow\Gamma_{z}$ are the injective homomorphisms given in the definition of orbifold. It is important to note that this condition is more than just an isomorphism of groups, but is an equality of their representations as actions on $\tilde V_z$. The collection of sets of this type form a subbasis for the corresponding $C^s$ topology on $C^r_{\COrb}(\orbify{O},\orbify{P})$.

Similarly, a $C^s$ neighborhood of $\lozengefunc f=(f,\{\Theta_{f,x}\})$ is defined to be
\begin{multline*}
\mathscr{N}^s(\lozengefunc f,\varepsilon)=\{\lozengefunc g=(g,\{\Theta_{g,x}\})\in\CROrbMaps^r(\orbify{O},\orbify{P})\mid\\%
 \redfunc g\in\mathscr{N}^s(\redfunc f,\varepsilon)\text{ and }%
\theta_{f(x)z}\circ\Theta_{f,x}=\theta_{g(x)z}\circ\Theta_{g,x}\text{ for all }x\in\orbify{O}\}
\end{multline*}
where
$\mathscr{N}^s(\redfunc f,\varepsilon)=\{\redfunc g\in\RedOrbMaps^r(\orbify{O},\orbify{P})\mid g\in\mathscr{N}^s(I\circ f,\varepsilon)%
\text{ for some }I\in\mathscr{ID}_{\orbify{N}}\}$.
Here, $\orbify{N}$ denotes, as usual, an open neighborhood of the image $f(\orbify{O})$.
\end{definition}

\begin{observation}\label{COrbMapNeighborhoodHomoObs} Suppose $\starfunc f=(f,\{\tilde f_x\},\{\Theta_{f,x}\})$ and 
$\starfunc f'=(f,\{\tilde f_x\},\{\Theta'_{f,x}\})$ are two complete orbifold maps in $C^r_{\COrb}(\orbify{O},\orbify{P})$ such that $f=f'$ as orbifold maps. Then
$\starfunc f'=(f,\{\tilde f_x\},\{\Theta'_{f,x}\})\notin\mathscr{N}^s(\starfunc f,\varepsilon)$ for any $\varepsilon$ unless $\{\Theta'_{f,x}\}=\{\Theta_{f,x}\}$. Otherwise, it would follow that
$\theta_{f(x)z}\circ\Theta_{f,x}=\theta_{f(x)z}\circ\Theta'_{f,x}$, contradicting injectivity of $\theta_{f(x)z}$. Of course the same argument shows that if $g=g'$ as orbifold maps, then $\starfunc g=(g,\{\tilde g_x\},\{\Theta_{g,x}\})$ and 
$\starfunc g'=(g,\{\tilde g_x\},\{\Theta'_{g,x}\})$ cannot both belong to a neighborhood 
$\mathscr{N}^s(\starfunc f,\varepsilon)$ unless $\{\Theta_{g,x}\}=\{\Theta'_{g,x}\}$. 
As a consequence, we see that 
the preimage $q^{-1}\left(\mathscr{N}^s(f,\varepsilon)\right)\subset C^r_{\COrb}(\orbify{O},\orbify{P})$ is a (finite) disjoint union of neighborhoods of the form
$\mathscr{N}^s(\starfunc f_i,\varepsilon)$ where $\starfunc f_i=(f,\{\tilde f_x\},\{\Theta_{f,x}\}_i)$. Similarly, we see that the preimage $q_{\ssstriangledown}^{-1}\left(\mathscr{N}^s(\redfunc f,\varepsilon)\right)\subset\CROrbMaps^r(\orbify{O},\orbify{P})$ is a (finite) disjoint union of neighborhoods of the form
$\mathscr{N}^s(\lozengefunc f_i,\varepsilon)$ where $\lozengefunc f_i=(f,\{\Theta_{f,x}\}_i)$.
\end{observation}

For reference, we have the following diagram of maps:

\begin{equation*}\label{LocalMapDiagram}
      \xymatrix{
     &  {\mathscr{N}^s(\starfunc f,\varepsilon)}\ar[dl]_{q_{\scriptscriptstyle\blacklozenge}}\ar[dr]^{q}& \\
     {\mathscr{N}^s(\lozengefunc f,\varepsilon)}\ar[rdd]_{q_{\scriptscriptstyle\blacktriangledown}}\ar[dr]^{q_\dagger}  & & %
     {\mathscr{N}^s(f,\varepsilon)}\ar[ldd]^{q_\bullet}\\
     & {\mathscr{N}^s([\lozengefunc f],\varepsilon)}\ar[d]_{q_\ddagger} & &\\
     & {\mathscr{N}^s(\redfunc f,\varepsilon)} & 
     }
\end{equation*}

We now show that the action of identity maps is compatible with the $C^s$ topology on $\COrbMaps^r(\orbify{O},\orbify{P})$. Let 
$\starfunc f=(f,\{\tilde f_x\},\{\Theta_{f,x}\})\in\COrbMaps^r(\orbify{O},\orbify{P})$ and let the orbifold $\orbify{N}$ be an open neighborhood of the image of $f(\orbify{O})$. Let $\varepsilon>0$ be chosen so that if 
$\starfunc g=(g,\{\tilde g_x\},\{\Theta_{g,x}\})\in\mathscr{N}^s(\starfunc f,\varepsilon)$, then $g(\orbify{O})\subset\orbify{N}$. Let $I\in\mathscr{ID}_\orbify{N}$ 
be an orbifold map lift of the identity over $\orbify{N}$. Then by example~\ref{IdentityMap}, $I$ has a representation as
$I=(\textup{Id},\{\tilde w\mapsto\gamma_z\cdot\tilde w\},\{\delta\mapsto\gamma_z\delta\gamma_z^{-1}\})$.

\begin{lemma}\label{IdCompatCrTopology} With $\starfunc g\in\mathscr{N}^s(\starfunc f,\varepsilon)$ as above, we have
$I\circ\starfunc g\in\mathscr{N}^s(I\circ \starfunc f,\varepsilon)$. Thus, the local action of $\mathscr{ID}_\orbify{N}$ on a neighborhood of $\starfunc f$ is continuous. 
In fact, $I:\mathscr{N}^s(\starfunc f,\varepsilon)\to\mathscr{N}^s(I\circ \starfunc f,\varepsilon)$,
$\starfunc g\mapsto I\circ\starfunc g$, is a homeomorphism. In fact, the map $I$ is a $C^\infty$ diffeomorphism.
\end{lemma}

\begin{proof} The homeomorphism claim is immediate from the definitions once one realizes that if $\tilde U_x$ is chosen as in definition~\ref{CrTopOnComplete}, so that $\theta_{f(x)z}\circ\Theta_{f,x}=\theta_{g(x)z}\circ\Theta_{g,x}$, then for $\delta\in\Gamma_x$,
\begin{align*}
I\circ\starfunc f &= (f,\{\gamma_{z}\cdot\tilde f_x\},\{\gamma_{z}(\theta_{f(x)z}\circ\Theta_{f,x}(\delta))\gamma_{z}^{-1}\})\text{ and,}\\
I\circ\starfunc g &= (g,\{\gamma_{z}\cdot\tilde g_x\},\{\gamma_{z}(\theta_{g(x)z}\circ\Theta_{g,x}(\delta))\gamma_{z}^{-1}\}).
\end{align*}
The smoothness claims follow from lemma~\ref{ActionOfIdIsSmooth} in section~\ref{InfiniteDimAnalysis}.
\end{proof}

Assuming theorem~\ref{MainTheorem}, we can now prove corollary~\ref{CompleteReducedOrbifoldMapStructure}.

\begin{proof}[Proof of corollary~\ref{CompleteReducedOrbifoldMapStructure}]
Let $\lozengefunc f=(f,\{\Theta_{f,x}\})\in \CROrbMaps^r(\orbify{O}, \orbify{P})$. For small enough $\varepsilon>0$, it is clear that 
$\mathscr{Q}_{\star}=q_{\ssslozenge}^{-1}\left(\mathscr{N}^r(\lozengefunc f,\varepsilon)\right)=%
\bigsqcup\mathscr{N}^r(\starfunc f_i,\varepsilon)\subset\COrbMaps^r(\orbify{O}, \orbify{P})$,
is a finite disjoint union of neighborhoods of the form 
$\mathscr{N}^r(\starfunc f_i,\varepsilon)$ where 
$\starfunc f_i=(f,\{\tilde f_x\}_i,\{\Theta_{f,x}\})$. If $C_z=C_{\Gamma_z}(\theta_{f(x)z}\circ\Theta_{f,x}(\Gamma_x))$ 
(which $=C_{\Gamma_z}(\theta_{g(x)z}\circ\Theta_{g,x}(\Gamma_x))$), then 
lemmas~\ref{PointwiseCOrbQuotientCROrb} and \ref{IdCompatCrTopology} imply that 
$(\mathscr{ID}_{\orbify{N}})_{\{C_z\}}$ acts smoothly and transitively on fibers of $\mathscr{Q}_{\star}$.
Example~\ref{PointwiseCOrbQuotientCROrbNotFree} shows that this action is not necessarily free. To understand what happens under these circumstances, suppose that $I\in(\mathscr{ID}_{\orbify{N}})_{\{C_z\}}$ fixes 
$\starfunc f=(f,\{\tilde f_x\},\{\Theta_{f,x}\})\in q_{\ssslozenge}^{-1}(\lozengefunc f)$. Let 
$\starfunc g=(g,\{\tilde g_x\},\{\Theta_{g,x}\})\in\mathscr{N}^r(\starfunc f,\varepsilon)$. By lemma~\ref{IdCompatCrTopology},
$I\circ\starfunc g\in\mathscr{N}^r(\starfunc f,\varepsilon)$. This shows that $q_{\ssslozenge}$ defines an orbifold chart.
Since each $\mathscr{N}^r(\starfunc f_i,\varepsilon)$ is an open manifold by theorem~\ref{MainTheorem} and the action of
$(\mathscr{ID}_{\orbify{N}})_{\{C_z\}}$ is smooth, corollary~\ref{CompleteReducedOrbifoldMapStructure} follows.
\end{proof}

\begin{proposition}\label{qDaggerIsLocalHomeo} The quotient map $q_\dagger:\mathscr{N}^s(\lozengefunc f,\varepsilon)\to\mathscr{N}^s([\lozengefunc f],\varepsilon)$, $\lozengefunc f=(f,\{\Theta_{f,x}\})\mapsto [\lozengefunc f]=(f,\{[\Theta_{f,x}]\})$ is a local homeomorphism. In fact, it is the quotient map defined by the group action of $\mathscr{ID}_{\orbify{N}}$ acting via $\lozengefunc f\mapsto I\circ\lozengefunc f$.
\end{proposition}
\begin{proof} It is clear from the definitions that $\mathscr{Q}_\dagger=q_\dagger^{-1}(\mathscr{N}^s([\lozengefunc f],\varepsilon))$ consists of finite disjoint union of neighborhoods of the form $\mathscr{N}^s(\lozengefunc f_i,\varepsilon)$ where $\lozengefunc f_i=(f,\{\eta_{x,i}\Theta_{f,x}\eta_{x,i}^{-1}\})$. The last statement follows by observing that $\mathscr{ID}_{\orbify{N}}$ acts transitively on $\mathscr{Q}_\dagger$ and 
if $\lozengefunc f = I\circ\lozengefunc f$, for $I=(\text{Id},\{\tilde y\to\eta_x\cdot\tilde y\})$, then 
$I\in (\mathscr{ID}_{\orbify{N}})_{\{C_z\}}$ where $C_z=C_{\Gamma_z}(\theta_{f(x)z}\circ\Theta_{f,x}(\Gamma_x))$. Since 
$C_z$ also is $=C_{\Gamma_z}(\theta_{g(x)z}\circ\Theta_{g,x}(\Gamma_x))$ for any $\lozengefunc g\in\mathscr{N}^s(\lozengefunc f,\varepsilon)$, we see that any such $I$ fixes pointwise the entire neighborhood $\mathscr{N}^s(\lozengefunc f,\varepsilon)$ and the result follows.
\end{proof}

Later we will have need to refer to the following useful fact about the relation between $\starfunc f$ and maps 
$\starfunc g\in\mathscr{N}^s(\starfunc f,\varepsilon)$:

\begin{lemma}\label{StrataImplicationsOfLocalNeighborhood} 
Let $\starfunc g=(g,\{\tilde g_x\},\{\Theta_{g,x}\})\in\mathscr{N}^s(\starfunc f,\varepsilon)$. Then for each $x\in\orbify{O}$, 
$\Theta_{f,x}=\theta_{g(x)f(x)}\circ\Theta_{g,x}:\Gamma_x\to\Gamma_{f(x)}$. Moreover, $\tilde f_x(\tilde x)$ and $\tilde g_x(\tilde x)$ both belong to the same connected (closed) stratum $\tilde V_{f(x)}^{\Theta_{f,x}(\Gamma_x)}=\{\tilde y\in\tilde V_{f(x)}\mid \delta\cdot\tilde y=\tilde y\text{ for all }\delta\in\Theta_{f,x}(\Gamma_x)\}$.
\end{lemma}
\begin{proof} In definition~\ref{CrTopOnComplete} we may choose $z=f(x)$. This yields the stated equality of homomorphisms immediately. Recall that
$\tilde\psi_{g(x)f(x)}$ denotes a lift of the inclusion map $\psi_{g(x)f(x)}:V_{g(x)}\hookrightarrow V_{f(x)}$ given in the definition of orbifold atlas. So, 
$\tilde\psi_{g(x)f(x)}\circ\tilde g_x:\tilde U_x\to\tilde g_x(\tilde U_x)\hookrightarrow\tilde V_{f(x)}$  is equivariant relative to 
$\theta_{g(x)f(x)}\circ\Theta_{g,x}$, which by hypothesis is the same as $\Theta_{f,x}$. Thus, for each $\gamma\in\Gamma_x$ we have
\begin{align*}
\tilde\psi_{g(x)f(x)}\circ\tilde g_x(\tilde x) &=\tilde\psi_{g(x)f(x)}\circ\tilde g_x(\gamma\cdot\tilde x)=%
\left(\theta_{g(x)f(x)}\circ\Theta_{g,x}(\gamma)\right)\cdot\tilde g_x(\tilde x)\\
&=\Theta_{f,x}(\gamma)\cdot\tilde g_x(\tilde x)
\end{align*}
from which it follows that $\tilde g_x(\tilde x)\in\tilde V_{f(x)}^{\Theta_{f,x}(\Gamma_x)}$.
\end{proof}

\section{Applications to the Orbifold Diffeomorphism Group}\label{DiffeosSection}

In this section, we show how the discussion of the previous sections applies to orbifold diffeomorphisms. For simplicity, we will continue to assume that the orbifold $\orbify{O}$ is compact. In \cite{BB2006}, we studied the group of orbifold diffeomorphisms $\OrbDiff^r(\orbify{O})$ and the reduced orbifold diffeomorphisms $\RedOrbDiff^r(\orbify{O})$ showing that each carried the structure of a (topological) Banach/\Frechet manifold. In fact we expressed $\RedOrbDiff^r(\orbify{O})$ as the quotient $\OrbDiff^r(\orbify{O})/\mathscr{ID}$, where, of course, 
$\mathscr{ID}\subset\OrbDiff^r(\orbify{O})$ represents the (finite) group of orbifold map lifts of the identity on $\orbify{O}$.

For diffeomorphism groups, it is not hard to see that
the group of complete orbifold diffeomorphisms $\COrbDiff^r(\orbify{O})$ may be regarded as the same as $\OrbDiff^r(\orbify{O})$ in much the same way that example~\ref{IdentityMap} illustrated the correspondence $\starfunc{\mathscr{ID}}\leftrightarrow\mathscr{ID}$. This follows from the proof of corollary~1.2 in \cite{BB2006}, where it is shown that if $f_1,f_2\in\OrbDiff^r(\orbify{O})$ represent the same reduced diffeomorphism $\redfunc f\in\RedOrbDiff^r(\orbify{O})$, then $f_1\circ f_2^{-1}\in\mathscr{ID}$. In the diffeomorphism case, one should note that since all homomorphisms $\Theta_{f,x}$ are actually isomorphisms and we assume isotropy groups act effectively, the behavior exhibited in lemmas~\ref{PointwiseCOrbQuotientOrb} and \ref{PointwiseConjCROrbQuotientRedOrb} cannot occur. There can never be multiple
$\Theta_{f,x}$'s corresponding to a particular local lift $\tilde f_x$. Collecting the results of example~\ref{IdentityMap} and lemmas~\ref{PointwiseCOrbQuotientCROrb}, \ref{PointwiseCROrbQuotientRedOrb} and \ref{IdCompatCrTopology}, and exploiting the fact that, in the case of diffeomorphism groups, we have a {\em global} ($C^\infty$-) smooth action of $\mathscr{ID}$, we get the following algebraic and topological structure result.

\begin{theorem}\label{DiffeoGroupTheorem} Let $\orbify{O}$ be a compact smooth $C^r$ orbifold. Then the following sequences are exact:
\begin{gather*}
1\longrightarrow\mathscr{ID}\longrightarrow\OrbDiff^r(\orbify{O})\longrightarrow\RedOrbDiff^r(\orbify{O})\longrightarrow 1\\
1\longrightarrow C(\mathscr{ID})\longrightarrow\OrbDiff^r(\orbify{O})\longrightarrow\CROrbDiff^r(\orbify{O})\longrightarrow 1\\
1\longrightarrow \lozengefunc{\mathscr{ID}}=\mathscr{ID}/C(\mathscr{ID})\longrightarrow\CROrbDiff^r(\orbify{O})\longrightarrow\RedOrbDiff^r(\orbify{O})\longrightarrow 1\\
\end{gather*}
where $C(\mathscr{ID})$ denotes the center of $\mathscr{ID}$. Moreover, each of the diffeomorphism groups $\OrbDiff^r(\orbify{O})$, 
$\CROrbDiff^r(\orbify{O})$ and $\RedOrbDiff^r(\orbify{O})$ carries the structure of a smooth $C^\infty$
Banach ($r<\infty$)/\Frechet ($r=\infty$) manifold. 
\end{theorem}

\begin{theorem}\label{DiffeoGroupsAreTopGroups} Each of the diffeomorphism groups $\OrbDiff^r(\orbify{O})$, 
$\CROrbDiff^r(\orbify{O})$ and $\RedOrbDiff^r(\orbify{O})$ is a topological group. That is, composition and inversion are continuous. Futhermore, when $r=\infty$, $\OrbDiff^\infty(\orbify{O})$, 
$\CROrbDiff^\infty(\orbify{O})$ and $\RedOrbDiff^\infty(\orbify{O})$ are convenient \Frechet Lie groups.
\end{theorem}

\begin{proof} For $0<r<\infty$, the group multiplication $\mu(f,g)=f\circ g$ and inversion $\text{inv}(f)=f^{-1}$ in the diffeomorphism group
$\OrbDiff^r(\orbify{O})$, corresponds to composition and inversion of the $C^r$ local equivariant lifts. These operations are known only to be $C^0$. This follows by the so-called $\Omega$-lemma of Palais \cite{MR0248880}, a suitable version of which is stated as lemma~\ref{OmegaLemma} for completeness. Thus, $\OrbDiff^r(\orbify{O})$ is a topological group. The structure result of theorem~\ref{DiffeoGroupTheorem} then yields the topological group structure for
$\CROrbDiff^r(\orbify{O})$ and $\RedOrbDiff^r(\orbify{O})$.

For $r=\infty$, by lemma~\ref{CompositionIsSmooth}, group multiplication $\mu(f,g)=f\circ g$ is smooth since 
$\OrbDiff^\infty(\orbify{O})$ is an open submanifold of $\OrbMaps^\infty(\orbify{O},\orbify{O})$ by \cite{BB2006}*{section~7}. To show that inversion is smooth, we use the argument given in \cite{MR1471480}*{Theorem~43.1}. Let $c=(c^t,\{\tilde c^t_x\}):\R\to\OrbDiff^\infty(\orbify{O})\subset\OrbMaps^\infty(\orbify{O},\orbify{O})$ be a smooth curve. Then, in a local orbifold chart by corollary~\ref{CurveSmoothEquivLiftSmooth}, the mapping 
$\tilde c_x^\wedge:(0,1)\times\tilde U_x\to\tilde V_z$ is smooth and $(\text{inv}\circ\tilde c_x)^\wedge$ satisfies the finite dimensional implicit equation $\tilde c_x^\wedge(t,(\text{inv}\circ\tilde c_x)^\wedge(t,\tilde y))=\tilde y$ for all $t\in\R$ and $\tilde y\in\tilde U_x$.
By the finite dimensional implicit function theorem, $(\text{inv}\circ\tilde c_x)^\wedge$ is smooth in $(t,\tilde y)$. Hence, by corollary~\ref{CurveSmoothEquivLiftSmooth}, $\text{inv}$ maps smooth curves to smooth curves and is thus smooth. This shows that $\OrbDiff^\infty(\orbify{O})$ is a convenient \Frechet Lie group and thus, by the structure results of theorem~\ref{DiffeoGroupTheorem}, so are
$\CROrbDiff^\infty(\orbify{O})$ and $\RedOrbDiff^\infty(\orbify{O})$.
\end{proof}

\section{Why Non-Orbifold Structure Stratifications Arise}\label{StratificationExampleSection}

In this section, we wish to give an example on why non-orbifold structure stratifications arise in the topological structure of our orbifold maps. We first recall a definition of stratification in the infinite-dimensional setting. We will use the definition found in 
\cite{MR0347323} or \cite{MR0418147} for infinite-dimensional stratifications although we do not need the full generality presented in these references. In our case, each point with a stratified neighborhood has only a finite number of strata coming together. 

\begin{definition}[\cites{MR0347323, MR0418147}]\label{StratificationDef} Let $X$ be a topological space and $\mathscr{A}$ a countable set with partial order $\prec$. A \emph{partition} of $X$ is a collection of non-empty pairwise disjoint subspaces $\{X_\alpha\}$ indexed by $\mathscr{A}$ such that 
$X=\cup_{\alpha\in\mathscr{A}}X_\alpha$. A partition $\{X_\alpha\}_{\alpha\in\mathscr{A}}$ is a \emph{stratification} of $X$ if 
\begin{enumerate}
\item\label{StrataCondition1} each $X_\alpha$ is a submanifold when given the topology induced by $X$ and,
\item\label{StrataCondition2} $X_\alpha\cap\overline{X_\beta}\ne\emptyset$, $\alpha\ne\beta$, then 
$\beta\prec\alpha$ and $X_\alpha\subset \overline{X_\beta}$.
\end{enumerate}
The $X_\alpha$ are called the \emph{strata} of the stratification and may have many connected components. Moreover, condition \eqref{StrataCondition2} implies that $\overline{X_\beta}-X_\beta\subset\cup_{\alpha\succ\beta}X_\alpha$.
\end{definition}

Before we show how these stratifications arise, we first present a simple example to help motivate the discussion.

\begin{example}\label{StratAriseEx1} Consider the situation described in Example~\ref{RZ2Example}:
$\orbify{O}$ is the orbifold $\R/\Z_2$ where $\Z_2$ acts on $\R$ via $x \to -x$ and $f:\orbify{O}\to\orbify{O}$ is the
constant map $f\equiv 0$. The map $\tilde f_0\equiv 0$ is a local equivariant lift of $f$ at
$x=0$ using either of the homomorphisms $\Theta_{f,0}=\textup{Id}$
or $\Theta'_{f,0}\equiv e$. Of course, for $x\ne 0$, we set $\tilde f_x\equiv 0$ and 
$\Theta_{f,x}=\Theta'_{f,x}=\text{ the trivial homomorphism }\Gamma_x=e\mapsto e\in\Gamma_0=\Z_2$. 
Thus, we have two complete orbifold maps $\starfunc f=(f,\{\tilde f_x\},\{\Theta_{f,x}\})$ and $\starfunc f'=(f,\{\tilde f_x\},\{\Theta'_{f,x}\})$ which cover the same orbifold map $f=(f,\{\tilde f_x\})$. 

We need to first compute $\mathscr{N}^r(\starfunc f,\varepsilon)$. We will do this in detail since this is the first time we have done an explicit computation of this type. Using definition~\ref{CrTopOnComplete} and the notation there, let 
$\starfunc g\in\mathscr{N}^r(\starfunc f,\varepsilon)$. For all $x\in\orbify{O}$ we may choose $z=0$ and thus $\tilde V_z=\tilde V_0$ may be chosen to be the interval $(-\varepsilon,\varepsilon)$ as a chart about $0$ in the target. There are two cases to consider: $x=0$ and $x\ne 0$. For $x=0$, let $\tilde U_0$ be any orbifold chart about $0$. It follows that the local lift $\tilde g_0$ over $x=0$ must take $0\in\tilde U_0$ to $0\in\tilde V_0$.
To see this, suppose to the contrary that $\tilde g_0(0)=\tilde y\ne 0$. By definition~\ref{CrTopOnComplete}, we must have the following equality of homomorphisms from $\Z_2=\Gamma_0$ to $\Gamma_0$:
\begin{align*}
\theta_{f(0)0}\circ\Theta_{f,0}& =\theta_{g(0)0}\circ\Theta_{g,0}\Longleftrightarrow\\
\textup{Id}& =\theta_{y0}\circ\Theta_{g,0}
\end{align*}
However, $\Theta_{g,0}:\Gamma_0\to\Gamma_y=\{e\}$ has nontrivial kernel which contradicts the last line above. We thus may conclude that for $x=0$,
$\tilde g_0(0)=0$ and $\Theta_{g,0}=\Theta_{f,0}=\text{Id}$. From $\Theta_{g,0}=\text{Id}$, it follows that the local lift
$\tilde g_0$ must be an odd function. For $x\ne 0$, there is no restriction on $\tilde g_x$ arising from equivariance since $\Gamma_x=\{e\}$ and
$\theta_{f(x)0}\circ\Theta_{f,x}=\theta_{g(x)0}\circ\Theta_{g,x}:\Gamma_x\to\Gamma_0$ will always be the trivial homomorphism $e\mapsto e$.
Putting this all together we have shown that
\begin{align*}
q\left(\mathscr{N}^r(\starfunc f,\varepsilon)\right)&=\{g\in\OrbMaps^r(\orbify{O})\mid \|\tilde{g}_x\|<\varepsilon\text{ and }\tilde{g}_0\text{ is an odd function}\}.
\end{align*}
We now use a similar argument to compute $\mathscr{N}^r(\starfunc f',\varepsilon)$. Let $\starfunc g'\in\mathscr{N}^r(\starfunc f',\varepsilon)$.
For $x\ne 0$, $\Theta'_{f,x}=\Theta_{f,x}$,  so we conclude as above that there is no restriction on $\tilde g'_x$ arising from equivariance. On the other hand, for $x=0$ we must have the equality of homomorphisms
$\theta_{f(0)0}\circ\Theta'_{f,0} =\theta_{g'(0)0}\circ\Theta_{g',0}:\Gamma_0\to\Gamma_0$. Since $\Theta'_{f,0}\equiv e$, injectivity of
$\theta_{g'(0)0}$ implies that $\Theta_{g',0}\equiv e$. Thus, there is no restriction on $\tilde g'_0$ arising from equivariance either and we can conclude that
$$q\left(\mathscr{N}^r(\starfunc f',\varepsilon)\right)=\{g'\in\OrbMaps^r(\orbify{O})\mid \|\tilde{g}'_x\|<\varepsilon\}.$$
Here, it is clear that $q\left(\mathscr{N}^r(\starfunc f,\varepsilon)\right)$ is a proper subset (later, a submanifold) of 
$q\left(\mathscr{N}^r(\starfunc f',\varepsilon)\right)$, and that any orbifold map $g\in\mathscr{N}^r(f,\varepsilon)$ must be in 
$q\left(\mathscr{N}^r(\starfunc f',\varepsilon)\right)$ so that the topological structure of a neighborhood of $f$ is completely determined from only an understanding of the topological structure of 
$q\left(\mathscr{N}^r(\starfunc f',\varepsilon)\right)$ which in turn is determined by the structure of $\mathscr{N}^r(\starfunc f',\varepsilon)$, which will be shown to be a manifold.
\end{example}

Unfortunately, in general, the topological structure of a neighborhood of an orbifold map $f$ is rarely determined completely by the topological structure of a {\em single} neighborhood of one of its complete orbifold lifts $q^{-1}(f)$. This is illustrated in the next example.

\begin{figure}[h]
   \centering
   \includegraphics[scale=.5]{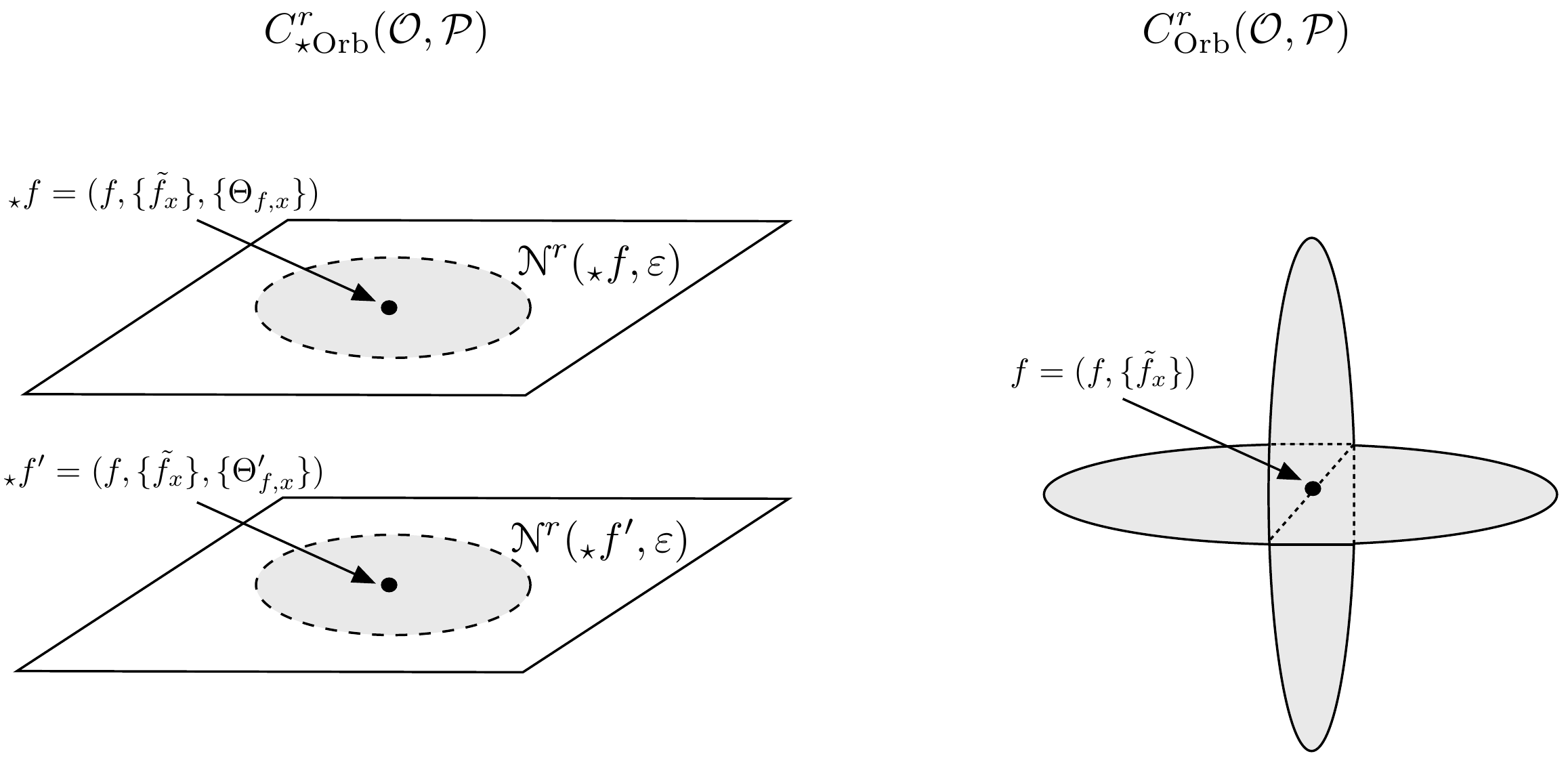} 
   \caption{A stratified neighborhood}
   \label{StratifiedNbhd}
\end{figure}

\begin{example}\label{StratAriseEx2}
Let $\orbify{O} = \R/\Z_2$ with $\Z_2$ acting with generator $\alpha$, where $\alpha\cdot x = -x$ as above. Let
$\orbify{P} = \R^3/(\Z_2\times \Z_2)$ where
$\Z_2\times \Z_2 = \langle j, k \mid j^2 = k^2 = 1 = \left[j,k\right]\rangle$ with the action defined by $j\cdot(x, y, z) = (-x, y, -z)$ and
$k\cdot(x, y, z) = (-x, -y, z)$.  Consider $\redfunc f\in\RedOrbMaps^r(\orbify{O},\orbify{P})$ defined by $\redfunc f(y_1) = (y_1, 0, 0)$ and choose
the orbifold map $f\in q_\bullet^{-1}(\redfunc f)$ given by $f=(f,\{\tilde y_1\mapsto (\tilde y_1,0,0)\})$. That is, 
for each $x\in\orbify{O}$, the local lift $\tilde f_x(\tilde y_1)=(\tilde y_1,0,0)$ on $\tilde U_x$. Since $\Gamma_x$ is trivial when $x\ne0$ and
$\Gamma_0=\Z_2$, there are precisely two complete maps in $q^{-1}(f)$:
\begin{eqnarray*}
 \starfunc{f} &=& (f, \{\tilde y_1 \mapsto (\tilde y_1, 0, 0)\}, \{\Theta_{f,0}:\alpha\mapsto j\})\\
 \starfunc{f}' &=& (f, \{\tilde y_1 \mapsto (\tilde y_1, 0, 0)\}, \{\Theta'_{f,0}:\alpha\mapsto k\}).
\end{eqnarray*}
Note that since $\Theta_{f,x}$ is the trivial homomorphism $e\mapsto e$ for all $x\ne 0$, we have only indicated the two possible homomorphisms at $x=0$,
namely, $\Theta_{f,0},\Theta'_{f,0}:\Gamma_0=\Z_2\to\Gamma_{(0,0,0)}=\Z_2\times \Z_2$.

We will proceed as in example~\ref{StratAriseEx1} and first compute $\mathscr{N}^r(\starfunc f,\varepsilon)$. Let 
$\starfunc g\in\mathscr{N}^r(\starfunc f,\varepsilon)$. Then $\starfunc g$ has a representation 
$$\starfunc g=(g,\{\tilde g_x=(\tilde y_1+(\tilde{g}_1)_x(\tilde y_1),(\tilde{g}_2)_x(\tilde y_1),(\tilde{g}_3)_x(\tilde y_1)\},\{\Theta_{g,x}\}).$$
For $x\ne 0$ we have $\Gamma_x=\{e\}$ so, like before, there is no restriction on $(\tilde{g}_i)_x, i=1,2,3$ arising from equivariance. Thus, we focus on lifts
$(\tilde{g}_i)_0$ over a chart $\tilde U_0$ about $x=0$. We may assume that $\tilde V_z=\tilde V_{\vec{0}}$ where we have shortened the subscript $(0,0,0)$ to $\vec{0}\in\R^3$. We will continue to do this for the remainder of this example. 

By lemma~\ref{StrataImplicationsOfLocalNeighborhood}, $\tilde g_0(0)\in\tilde V_{\vec{0}}^{\{e,j\}}=y$-axis and $\Theta_{g,0}:\Gamma_0\to\Gamma_{g(0)}$ is $\alpha\mapsto j$. We now compute
\begin{align*}
\tilde g_0(\alpha\cdot\tilde y_1)&=\tilde g_0(-y_1)=%
\left(-\tilde y_1+(\tilde{g}_1)_0(-\tilde y_1),(\tilde{g}_2)_0(-\tilde y_1),(\tilde{g}_3)_0(-\tilde y_1)\right).\\
\intertext{On the other hand,}
\tilde g_0(\alpha\cdot\tilde y_1) &=\Theta_{g,0}(\alpha)\cdot\tilde g_0(y_1)= j\cdot\tilde g_0(y_1)\\
& = \left(-\tilde y_1-(\tilde{g}_1)_0(\tilde y_1),(\tilde{g}_2)_0(\tilde y_1),-(\tilde{g}_3)_0(\tilde y_1)\right).
\end{align*}
Thus,
\begin{multline*}q(\mathscr{N}^r(\starfunc f,\varepsilon))=\{g\in\OrbMaps^r(\orbify{O},\orbify{P})\mid \\
\|\tilde g_x-\tilde f_x\|<\varepsilon\text{ with } (\tilde{g}_1)_0,(\tilde{g}_3)_0\text{ odd functions and }%
(\tilde{g}_2)_0\text{ an even function}\}.
\end{multline*}
Similarly, we have
\begin{multline*}q(\mathscr{N}^r(\starfunc f',\varepsilon))=\{g'\in\OrbMaps^r(\orbify{O},\orbify{P})\mid \\
\|\tilde g'_x-\tilde f_x\|<\varepsilon\text{ with } (\tilde{g}'_1)_0,(\tilde{g}'_2)_0\text{ odd functions and }%
(\tilde{g}'_3)_0\text{ an even function}\}.
\end{multline*}
Thus, the corresponding neighborhood of the orbifold map $f$ is the union of two sets
$\mathscr{N}^r(f,\varepsilon)=q(\mathscr{N}^r(\starfunc f,\varepsilon))\cup q(\mathscr{N}^r(\starfunc f',\varepsilon))$ each of which will later be shown to carry a Banach/\Frechet manifold structure. Their intersection is along the submanifold
\begin{multline*}
\mathcal{H}=q(\mathscr{N}^r(\starfunc f,\varepsilon))\cap q(\mathscr{N}^r(\starfunc f',\varepsilon))=\\
\{h\in\OrbMaps^r(\orbify{O},\orbify{P})\mid\|\tilde h_x-\tilde f_x\|<\varepsilon
\text{ with } \tilde h_0(\tilde y_1)=(\tilde y_1+(\tilde{h}_1)_0(\tilde y_1),0,0)\\\text{ where }(\tilde{h}_1)_0\text{ is an odd function}\}.
\end{multline*}
Thus, the neighborhood $\mathscr{N}^r(f,\varepsilon)$ has a stratified structure (see figure~\ref{StratifiedNbhd}): Just let $\mathscr{A}=\{\alpha,\beta,\gamma\}$ with partial order 
$\beta\prec\alpha$, $\gamma\prec\alpha$ and define $X=\mathscr{N}^r(f,\varepsilon)$, $X_\alpha=\mathcal{H}$, 
$X_\beta=q(\mathscr{N}^r(\starfunc f,\varepsilon))-\mathcal{H}$, and $X_\gamma=q(\mathscr{N}^r(\starfunc f',\varepsilon))-\mathcal{H}$.
Moreover, since $\mathscr{N}^r(f,\varepsilon)-\orbify{H}$ is not connected we see that this stratified structure is not that of an orbifold structure as removal of the singular set of an orbifold never disconnects a connected component of the orbifold \cites{BorPhD, MR1218706}. Furthermore, if we let $\orbify{N}$ denote an open neighborhood of the image $f(\orbify{O})$, then from \cite{BB2006} a neighborhood of the reduced orbfold map $\redfunc f$ is given by 
$\mathscr{N}^r(\redfunc f,\varepsilon)=\mathscr{N}^r(f,\varepsilon)/\mathscr{ID}_{\orbify{N}}$ where $\mathscr{ID}_{\orbify{N}}$ acts in such a way that the quotient map restricts on each stratum to give a smooth orbifold chart (see proof of corollary~\ref{ReducedOrbifoldMapStructure} which appears at the end of
section~\ref{StratifiedNeighborhoodsSection}). Thus,
$\mathscr{N}^r(\redfunc f,\varepsilon)$ has a non-orbifold structure stratification also.
\end{example}

\section{The Tangent Orbibundle, Pullbacks and Orbisections}\label{TangentOrbibundleSection}

\subsection*{The tangent orbibundle} We recall the definition of the tangent orbibundle of a smooth $C^{r+1}$ orbifold.

\begin{definition}\label{tangentorbibundle}
  Let $\orbify{O}$ be an $n$--dimensional $C^{r+1}$ orbifold.  The
  \emph{tangent orbibundle} of $\orbify{O}$,
  $p:T\orbify{O}\to\orbify{O}$, is the $C^{r}$ orbibundle defined as
  follows.  If $(\tilde U_x,\Gamma_x)$ is an orbifold chart around
  $x\in \orbify{O}$ then $p^{-1}(U_x)\cong(\tilde
  U_x\times\R^n)/\Gamma_x$ where $\Gamma_x$ acts on $\tilde
  U_x\times\R^n$ via $\gamma\cdot(\tilde y, \tilde v) =
  \left(\gamma\cdot\tilde y, d\gamma_{\tilde y}(\tilde v)\right)$. In
  keeping with tradition, we denote the fiber $p^{-1}(x)$ over $x\in
  U_x$ by $T_x\orbify{O}\cong\R^n/\Gamma_x$. Note that, in general, if
  $\Gamma_x$ is non-trivial then $T_x\orbify{O}$ will be a convex cone
  rather than a vector space. Locally we have the diagram:
  
  \begin{equation*}
    \xymatrix{
      {T\tilde U_x\cong\tilde U_x\times\R^n}\ar[rrr]^{\Pi_x}\ar[d]_{\text{pr}_1} & & &{(\tilde U_x\times\R^n)/\Gamma_x}\ar[d]^{p}\\
      {\tilde U_x}\ar[rrr]^{\pi_x}& &  & {U_x}
    }  \end{equation*}
  where $\text{pr}_1:\tilde U_x\times\R^n\to\tilde U_x$ denotes the
  projection onto the first factor $(\tilde y,\tilde v)\mapsto\tilde
  y$ (which is a specific choice of lift of $p$).
\end{definition}

\subsection*{Pulling back an orbibundle}
The definition of the pullback of an orbibundle depends crucially on
the notion of orbifold map.  In simple examples, we will see that a unique notion of pullback exists only when using
complete orbifold maps. On the other hand, we will see that once one has a pullback bundle defined via a complete orbifold
map $\starfunc{f}$, there is no difference between the notion of
an orbisection and a complete orbisection.  Not surprisingly, if one tries to define a useful notion of reduced or complete reduced
orbisection one loses the vector space structure on the
space of such sections. As in the case of the tangent orbibundle, the pullback bundle will be an example of the more general notion of a linear orbibundle given in \cite{MR1926425}.

\begin{definition}\label{PullbackBundleDef}
 Let $\orbify{O}, \orbify{P}$ be $C^{r+1}$ orbifolds of dimension $n$ and $m$, respectively. Given $\starfunc{f}\in \COrbMaps^{r+1}(\orbify{O}, \orbify{P})$ we
 define the pullback of the tangent orbibundle to $\orbify{P}$ by
 $\starfunc{f}$, $\starfunc{f}^*(T\orbify{P})$ as follows:  Let
 $\starfunc{f} = (f, \{\tilde f_x \}, \{\Theta_{f,x}\})$ and let $\tilde U_x$ and $\tilde V_{f(x)}$ be orbifold
 charts about $x\in\orbify{O}$ and
 $f(x)\in\orbify{P}$ respectively. Define the pullback $\starfunc{f}^*(T\orbify{P})$ to be the
 orbifold with charts of the form (a fibered product)
 \begin{equation*}
   \tilde U_x\times_{\tilde V_{f(x)}}T\tilde V_{f(x)} = %
   \{(\tilde y, \tilde\xi) \in \tilde U_x\times T\tilde V_{f(x)} \mid \tilde f_x(\tilde y) = \text{pr}_1(\tilde\xi)\}
 \end{equation*}
 where $\tilde f_x:\tilde U_x\to \tilde V_{f(x)}$ and $\text{pr}_1:T\tilde V_{f(x)}\to \tilde V_{f(x)}$ is the tangent bundle
 projection.  If we write $\tilde\xi=[\tilde f_x(\tilde y),\tilde v]\in T\tilde V_{f(x)}=\tilde V_{f(x)}\times\R^m$, the action of $\Gamma_x$ is specified in local
 coordinates by:
 \begin{align*}
   \gamma\cdot(\tilde y, \tilde\xi) & = (\gamma\cdot\tilde y, \Theta_{f,x}(\gamma)\cdot\tilde\xi)=(\gamma\cdot\tilde y,[\Theta_{f,x}(\gamma)\cdot\tilde f_x(\tilde y), d(\Theta_{f,x}(\gamma))_{\tilde f_x(\tilde y)}\cdot\tilde v])\\
   &=(\gamma\cdot\tilde y,[\tilde f_x(\gamma\cdot \tilde y), d(\Theta_{f,x}(\gamma))_{\tilde f_x(\tilde y)}\cdot\tilde v])%
   \in\tilde U_x\times_{\tilde V_{f(x)}}T\tilde V_{f(x)}
 \end{align*}
 where $\tilde y\in \tilde U_x$ and $\tilde v\in \text{pr}_1^{-1}(\tilde f_x(\tilde y))$. Also, we let
 $\text{pr}_2:T\tilde V_{f(x)}\cong\tilde V_{f(x)}\times\R^m\to\R^m$, $\tilde\xi\mapsto \tilde v$ be the fiber projection. This gives 
 $\starfunc{f}^*(T\orbify{P})$ the structure of a smooth $C^r$ $m$-dimensional linear orbibundle over $\orbify{O}$. In an abuse of notation,
 $p:\starfunc{f}^*(T\orbify{P})\to \orbify{O}$ will denote the orbibundle projection. Denote the fiber over $x$, by 
 $p^{-1}(x)=\starfunc{f}^*(T\orbify{P})_x$. In local coordinates, we have the diagram (all vertical arrows are quotient maps by respective group actions):
 
 \begin{equation*}\label{PullbackOrbibundleDiagram}
  \xymatrix{
   & \tilde U_x\times_{\tilde V_{f(x)}}T\tilde V_{f(x)}\ar[dl]_{\textrm{pr}_1}\ar'[d][ddd]^(.3){\Pi_x} \ar[rrr]^{\tilde F_x} &  &  &%
   *{\parbox{1.5in}{{$\begin{aligned}\hspace{2ex}&T\tilde V_{f(x)}\cong\\(&\tilde V_{f(x)}\times\R^m)\end{aligned}$}}}%
   \ar[dl]_(.6){\textrm{pr}_1}\ar[dd]_{\Pi_{\Theta_{f,x}}}\ar@/^3.5pc/[ddd]_(.45){\Pi_{f(x)}}\\
    {\tilde U_x}\ar[rrr]^(.6){\tilde f_x}\ar[ddd]^{\pi_x}& & &{\tilde V_{f(x)}}\ar[ddd]^(.55){\pi_{f(x)}}\\
     & &  &  & T\tilde V_{f(x)}/\Theta_{f,x}(\Gamma_x)\ar[d]\\
     &*{\parbox{1.5in}{{$\begin{aligned}&\starfunc f^*(TV_{f(x)})\cong\\(\tilde U_x &\times_{\tilde V_{f(x)}}T\tilde V_{f(x)})/\Gamma_x\end{aligned}$}}}%
     \ar[dl]^p\ar@{-->}[urrr]^(.4){F} |!{[rr];[rrd]}\hole%
      &  & &*{\parbox{1.5in}{{$\begin{aligned}&TV_{f(x)}\cong\\(\tilde V_{f(x)}&\times\R^m)/\Gamma_{f(x)}\end{aligned}$}}}\ar[dl]^p\\
    U_x\ar[rrr]^{f} &  &  & V_{f(x)}\\
  }
\end{equation*}
The map $\tilde F_x:\tilde U_x\times_{\tilde V_{f(x)}}T\tilde V_{f(x)}\to T\tilde V_{f(x)}$ given by $\tilde F_x(\tilde y,\tilde\xi)=\tilde\xi$ induces a map $F:\starfunc f^*(TV_{f(x)})\to T\tilde V_{f(x)}/\Theta_{f,x}(\Gamma_x)$ defined by $F(y,\xi)=\Pi_{\Theta_{f,x}}\circ\tilde F_x\circ\Pi_x^{-1}(y,\xi)$. This is well defined since, for any $\gamma\in\Gamma_x$,
$\tilde F_x(\gamma\cdot (\tilde y,\tilde\xi))=\tilde F_x(\gamma\cdot\tilde y,\Theta_{f,x}(\gamma)\cdot\tilde\xi)=\Theta_{f,x}(\gamma)\cdot\tilde\xi$ and
$\Pi_{\Theta_{f,x}}(\Theta_{f,x}(\gamma)\cdot\tilde\xi)=\Pi_{\Theta_{f,x}}(\tilde \xi)$.
\end{definition}

Note that the pullback is defined only if we have all the
information contained in both the choices of local lifts $\{\tilde
f_x\}$ and the choices of the homomorphisms $\Theta_{f,x}\in\operatorname{Hom}(\Gamma_x, \Gamma_{f(x)})$.  That is, all of the information of a complete orbifold map is used. As an illustration of
the necessity for needing to use complete orbifold maps to define pullbacks we give two examples. The first example shows that, unless complete orbifold maps are used, the pullback orbibundle is not well-defined even up to a reasonable notion of equivalence. 

\begin{example}\label{RZ2ExamplePullback}
Consider the situation from example~\ref{StratAriseEx1}:
$\orbify{O}=\orbify{P}$ is the orbifold $\R/\Z_2$ where $\Z_2$ acts on $\R$ via $x \to -x$ and $f:\orbify{O}\to\orbify{P}$ is the
constant map $f\equiv 0$.  Note
that $T\orbify{O} = T\orbify{P} = \R^2/\Z_2$ where the generator $\alpha$
of $\Z_2$ acts via $(x, y)\mapsto (-x, -y)$ and the bundle projection
$p$ is just projection onto the first factor.  We note that for
$x\ne 0$, $p^{-1}(x) = \R$ and that $p^{-1}(0) = \R/\Z_2$.
Let $\starfunc f$ and $\starfunc f'$ be the two complete orbifold maps from 
example~\ref{StratAriseEx1} which cover the orbifold map $f=(f,\{\tilde f_x\})$. Then we claim that
 \begin{align*}
   \starfunc{f}^*(T\orbify{P}) &\cong T\orbify{O}\\
   \intertext{ while }
   \starfunc{f'}^*(T\orbify{P}) &\cong \orbify{O}\times \R
 \end{align*}
To see this we work in local coordinates: Since $\tilde f_x\equiv 0$, we may take $\tilde V_0$ as a chart about $f(x)$ for all $x$. Thus for each $x$,
\begin{align*}
\tilde U_x\times_{\tilde V_0}T\tilde V_0&=\{(\tilde y, \tilde\xi) \in \tilde U_x\times T\tilde V_0\mid 0=\text{pr}_1(\tilde \xi)\}\\
&=\tilde U_x\times T_0\tilde V_0\cong\tilde U_x\times\R
\end{align*}
Now for $x\ne 0$, $\Gamma_x=\{e\}$ and so the action of $\Gamma_x$ on $\tilde U_x\times_{\tilde V_0}T\tilde V_0$ is necessarily trivial. If we denote the orbibundle projections $p:\starfunc{f}^*(T\orbify{P})\to \orbify{O}$ and $p':\starfunc{f'}^*(T\orbify{P})\to \orbify{O}$, then
$p^{-1}(U_x)=p'^{-1}(U_x)\cong\tilde U_x\times\R$. On the other hand, for $x=0$, since $\Theta_{f,0}(\alpha)=\alpha$ and $\Theta'_{f,0}(\alpha)=e$ we see that
\begin{multline*}
p^{-1}(U_0)\cong(\tilde U_0\times\R)/\Gamma_0 \text{ where the action of $\Gamma_0$ is }\\
\alpha\cdot(\tilde y,\tilde v_0)=(\alpha\cdot\tilde y,d(\Theta_{f,0}(\alpha))_{\tilde f_0(\tilde y)}\cdot\tilde v_0)=
(\alpha\cdot\tilde y,d(-\text{Id})_{0}\cdot\tilde v_0)=(-\tilde y,-\tilde v_0)
\end{multline*}
while
\begin{multline*}
p'^{-1}(U_0)\cong(\tilde U_0\times\R)/\Gamma_0 \text{ where the action of $\Gamma_0$ is }\\
\alpha\cdot(\tilde y,\tilde v_0)=(\alpha\cdot\tilde y,d(\Theta'_{f,0}(\alpha))_{\tilde f_0(\tilde y)}\cdot\tilde v_0)=%
(\alpha\cdot\tilde y,d(\text{Id})_{0}\cdot\tilde v_0)=(-\tilde y,\tilde v_0)
\end{multline*}
which is enough to substantiate our claim. Note that these orbibundles are not equivalent in any reasonable sense.

\end{example}

To further illustrate the complexity involved in pulling back the
tangent bundle by an orbifold map, the following is instructive.

\begin{example}
 \label{pullback:R3modZ2timesZ2:to:R1modZ2}
Consider the situation from example~\ref{StratAriseEx2}: $\orbify{O} = \R/\Z_2$ with $\Z_2$ acting with generator $\alpha$, where $\alpha\cdot x = -x$ as above and $\orbify{P} = \R^3/(\Z_2\times \Z_2)$ where
$\Z_2\times \Z_2 = \langle j, k \mid j^2 = k^2 = 1 = \left[j,k\right]\rangle$ with the action defined by $j\cdot(x, y, z) = (-x, y, -z)$ and
$k\cdot(x, y, z) = (-x, -y, z)$. We consider the two complete orbifold maps $\starfunc f$ and $\starfunc f'$ from example~\ref{StratAriseEx2} which cover the orbifold map $f=(f,\{\tilde y_1\mapsto (\tilde y_1,0,0)\})$ where $\redfunc f(y_1) = (y_1, 0, 0)$. We have for all $x$, 
$\tilde U_x\times_{\tilde V_{f(x)}}T\tilde V_{f(x)}\cong\tilde U_x\times\R^3$. Like in example~\ref{RZ2ExamplePullback},
for $x\ne 0$, $\Gamma_x=\{e\}$, and so the action of $\Gamma_x$ on $\tilde U_x\times_{\tilde V_{f(x)}}T\tilde V_{f(x)}$ is necessarily trivial. If we denote, as before, the orbibundle projections $p:\starfunc{f}^*(T\orbify{P})\to \orbify{O}$ and $p':\starfunc{f'}^*(T\orbify{P})\to \orbify{O}$, then
$p^{-1}(U_x)=p'^{-1}(U_x)\cong\tilde U_x\times\R^3$.
On the other hand, for $x=0$, since $\Theta_{f,0}(\alpha)=j$ and $\Theta'_{f,0}(\alpha)=k$ we see that
\begin{multline*}
p^{-1}(U_0)\cong(\tilde U_0\times\R^3)/\Gamma_0 \text{ where the action of $\Gamma_0$ is }\\
\alpha\cdot(\tilde y,\tilde v_0)=(\alpha\cdot\tilde y,dj_{\tilde f_0(\tilde y)}\cdot\tilde v_0)=
(-\tilde y,-(\tilde v_1)_0,(\tilde v_2)_0,-(\tilde v_3)_0)
\end{multline*}
where $\tilde v_0=((\tilde v_1)_0,(\tilde v_2)_0,(\tilde v_3)_0)\in\R^3$. Similarly,
\begin{multline*}
p'^{-1}(U_0)\cong(\tilde U_0\times\R^3)/\Gamma_0 \text{ where the action of $\Gamma_0$ is }\\
\alpha\cdot(\tilde y,\tilde v_0)=(\alpha\cdot\tilde y,dk_{\tilde f_0(\tilde y)}\cdot\tilde v_0)=
(-\tilde y,-(\tilde v_1)_0,-(\tilde v_2)_0,(\tilde v_3)_0).
\end{multline*}
Although the pullback orbibundles $\starfunc{f}^*(T\orbify{P})$ and $\starfunc{f'}^*(T\orbify{P})$ are naturally isomorphic, we will later see that neighborhoods of the zero section are taken by the Riemannian exponential map to the neighborhoods $\mathscr{N}^r(\starfunc f,\varepsilon)$ and
$\mathscr{N}^r(\starfunc f',\varepsilon)$ of example~\ref{StratAriseEx2}, respectively. This illustrates why it is necessary to use \emph{complete} orbifold maps in order to fully understand the topological structure of a neighborhood $\mathscr{N}^r(f,\varepsilon)$ of an orbifold map.
\end{example}

\subsection*{Orbisections} We now define a natural notion of section of a linear orbibundle. For a definition, see for example \cite{MR1926425}.
\begin{definition}\label{OrbiSection}
  A $C^r$ \emph{orbisection} of a $m$-dimensional linear orbibundle $p:\orbify{E}\to\orbify{O}$
  is a $C^r$ orbifold map $\sigma:\orbify{O}\to \orbify{E}$ such that
  $p\circ \sigma = \text{Id}_\orbify{O}$ and for any $x\in\orbify{O}$
  and chart $U_x$ about $x$, we have $\text{pr}_1\circ\tilde\sigma_x =
  \text{Id}_{\tilde U_x}$. That is, we take the {\em identity} lift of the identity $\text{Id}_\orbify{O}$ in $\tilde U_x$.
  Locally we have the diagram:

\begin{equation*}
  \xymatrix{
   & & & & {\tilde U_x\times\R^m}\ar[dl]^{\textrm{pr}_1}\ar[dd]^{\Pi_x}\\
    {\tilde U_x}\ar[rrr]_{\textrm{Id}_{\tilde U_x}}\ar[dd]_{\pi_x}\ar[urrrr]^(.4){\tilde\sigma_x}& & &{\tilde U_x}\ar[dd]^{\pi_x}\\
     &  &  & &*{\parbox{1in}{{$\begin{aligned}p^{-1}(U_x)\cong\\[-.05in](\tilde U_x\times\R^m)/\Gamma_x\end{aligned}$}}}\ar[dl]^p\\
    U_x\ar[rrr]_{\textrm{Id}}\ar[urrrr]^(.4){\sigma} |!{[urrr];[rrr]}\hole&  &  & U_x\\
  }
\end{equation*}
\end{definition}

Note that the action of $\Gamma_x$ on $\tilde U_x\times\R^m$ is given as part of the data defining $\orbify{E}$.
Although, in general, the class of complete orbifold maps is different from the
class of orbifold maps, as in the case for diffeomorphisms (section~\ref{DiffeosSection}), in the case of orbisections of the pullback  of a tangent orbibundle, these notions coincide.

\begin{proposition}\label{COrbisectionIsOrbisection}
Let $\starfunc f:\orbify{O}\to\orbify{P}$ be a $C^{r+1}$ complete orbifold map between $C^{r+1}$ orbifolds and let $\starfunc f^*(T\orbify{P})$ denote the pullback of the tangent orbibundle. Let
$\sigma=(\sigma, \{\tilde \sigma_x\})$ be a $C^r$ orbisection of $\starfunc f^*(T\orbify{P})$.
Then there is a unique homomorphism $\Theta_{\sigma,x}$ for which $\tilde\sigma_x$ is
$\Gamma_x$ equivariant. In other words, the set of orbisections can be identified with the set of complete orbisections $\sigma\leftrightarrow\starfunc\sigma$.
\end{proposition}

\begin{proof}
Given an orbifold chart $\tilde U_x$ around $x$ and an orbibundle
chart for $\starfunc{f}^*(T\orbify{P})$ with local product
coordinates $(\tilde y, \tilde \xi)=(\tilde y,[\tilde f_x(\tilde y),\tilde v])\in\tilde U_x\times_{\tilde V_{f(x)}}T\tilde V_{f(x)}$, the
local action of $\Gamma_x$ in these coordinates is given by 
$\gamma\cdot(\tilde y, \tilde \xi) = (\gamma\cdot\tilde y, \Theta_{f,x}(\gamma)\cdot\tilde\xi)=(\gamma\cdot\tilde y,[\Theta_{f,x}(\gamma)\cdot\tilde f_x(\tilde y), d(\Theta_{f,x}(\gamma))_{\tilde f_x(\tilde y)}\cdot\tilde v])$.  
With respect to these local coordinates, $\tilde\sigma_x$ has the form 
$\tilde\sigma_x(\tilde y) = (\tilde y, [\tilde f_x(\tilde y),\tilde s_x(\tilde y)])$ and if
$\Theta_{\sigma, x}:\Gamma_x\to\Gamma_{\sigma(x)}=\Gamma_x$ is some homomorphism for which $\tilde \sigma_x$ is
equivariant with respect to, then
\begin{align*}
   \tilde \sigma_x(\gamma\cdot\tilde y) &= (\gamma\cdot \tilde y, [\tilde f_x(\gamma\cdot\tilde y),\tilde s_x(\gamma\cdot\tilde y)])&\\
   &= \Theta_{\sigma, x}(\gamma)\cdot(\tilde y, [\tilde f_x(\tilde y),\tilde s_x(\tilde y)]) &\\
   &=(\Theta_{\sigma, x}(\gamma)\cdot\tilde y, [(\Theta_{f,x}(\Theta_{\sigma, x}(\gamma))\cdot\tilde f_x(\tilde y),%
   d\left(\Theta_{f,x}(\Theta_{\sigma, x}(\gamma))\right)_{\tilde f_x(\tilde y)}\cdot\tilde s_x(\tilde y)]) &
\end{align*}
Therefore, since $\Gamma_x$ acts effectively on $\tilde U_x$, $\Theta_{\sigma, x}(\gamma) = \gamma$ and $\Theta_{\sigma, x}=\text{Id}:\Gamma_x\to\Gamma_x$ for all $\gamma\in\Gamma_x$ and $x\in\orbify{O}$.
Furthermore, we get the equivariance relation $\tilde s_x(\gamma\cdot\tilde y)=d(\Theta_{f,x}(\gamma))_{\tilde f_x(\tilde y)}(\tilde s_x(\tilde y))$.
\end{proof}

Just as in the case of orbisections of the tangent orbibundle, the set of orbisections of the pullback tangent orbibundle carry a vector space structure.

\begin{proposition}\label{OrbiSecsAreVecSpace}
  Let $\starfunc f\in\COrbMaps^{r+1}(\orbify{O},\orbify{P})$. The set $\mathscr{D}^r_{\Orb}(\starfunc f^*(T\orbify{P}))$ of $C^r$ orbisections of
  the the pullback tangent orbibundle $\starfunc f^*(T\orbify{P})$ is naturally a real vector
  space with the vector space operations being defined pointwise.
\end{proposition}

\begin{proof} The argument here is basically the same as the corresponding argument for orbisections of the tangent orbibundle \cite{BB2006}.
Let $\sigma\in\mathscr{D}^r_{\Orb}(\starfunc f^*(T\orbify{P}))$.
Let $\tilde\sigma_x$ be the lift of $\sigma$.  Then
$\tilde\sigma_x(\tilde y)=(\tilde y,[\tilde f_x(\tilde y),\tilde s_x(\tilde y)])$. By proposition~\ref{COrbisectionIsOrbisection},
$\tilde s_x:\tilde U_x\to\R^m$ is such that $\tilde s_x(\gamma\cdot\tilde y)=d(\Theta_{f,x}(\gamma))_{\tilde f_x(\tilde y)}(\tilde s_x(\tilde y))$. 
In particular, since $\tilde x$ is a fixed point of the $\Gamma_x$ action on $\tilde U_x$,
we have $\tilde s_x(\tilde x)=\tilde s_x(\gamma\cdot\tilde x)=%
d(\Theta_{f,x}(\gamma))_{\tilde f_x(\tilde x)}(\tilde s_x(\tilde x))$.  Thus $\tilde s_x(\tilde x)$ is
a fixed point of the (linear) action of $\Gamma_x$ on $\R^m$ as defined in the pullback orbibundle. Note
that the set of such fixed points forms a {\em vector subspace} of
$\R^m$. As a result we may define a real vector space structure on
$\mathscr{D}^r_{\Orb}(\starfunc f^*(T\orbify{P}))$ as follows: For
$\sigma_i\in\mathscr{D}^r_{\Orb}(\starfunc f^*(T\orbify{P}))$, let
$\tilde\sigma_{i,x}$ be local lifts at $x$ as above. Define
\begin{align*}
  (\sigma_1+\sigma_2)(y) & =
  \Pi_x\big((\tilde\sigma_{1,x}+\tilde\sigma_{2,x})(\tilde y)\big)=%
  \Pi_x\big((\tilde y,\tilde s_{1}({\tilde y})+\tilde s_{2}({\tilde y}))\big)=\sigma_1(y)+\sigma_2(y)\\
  (\lambda\sigma)(y) & = \Pi_x\big((\lambda\tilde\sigma_{x})(\tilde
  y)\big)=%
  \Pi_x\big((\tilde y,\lambda\tilde s({\tilde
    y}))\big)=\lambda(\sigma(y))
\end{align*}
\end{proof}

\begin{proposition}\label{OrbisectionsOfPullBackIsBanachSpace} Let $\starfunc f\in\COrbMaps^{r+1}(\orbify{O},\orbify{P})$ with  $\orbify{O}$ compact (without boundary). 
The inclusion map
$\mathscr{D}^r_{\Orb}(\starfunc f^*(T\orbify{P}))\hookrightarrow
C^r_{\Orb}(\orbify{O}, \starfunc f^*(T\orbify{P}))$ induces a separable Banach
space structure on $\mathscr{D}^r_{\Orb}(\starfunc f^*(T\orbify{P}))$ for $1\le r<\infty$ and a separable \Frechet space structure if $r=\infty$.
\end{proposition}
\begin{proof}
\begin{sloppypar}
The argument here is also similar to the corresponding argument for orbisections of the tangent orbibundle \cite{BB2006}.
Let $\mathcal{D} = \{D_i\}_{i = 1}^N$ be a cover of $f(\orbify{O})$ by
a finite number of compact orbifold charts over each of which the tangent orbibundle $T\orbify{P}$ is trivialized. Then the collection
$\mathcal{C}=\{C_i=f^{-1}(D_i)\}$ is a finite cover of $\orbify{O}$ by compact subsets. By reindexing and shrinking $D_i$ if necessary, we may assume each $C_i$ is connected and is contained in a orbifold chart of $\orbify{O}$ and so, $C_i\cong\tilde C_i/\Gamma_i$. Let $\tilde C_i\times_{\tilde D_i}T\tilde D_i\cong\tilde C_i\times\R^m$ be the corresponding
orbifold charts for $\starfunc f^*(T\orbify{P})$ with action of 
$\Gamma_i: \gamma\cdot(\tilde y,\tilde\xi)=(\gamma\cdot\tilde y,%
[\Theta_{f,i}(\gamma)\cdot\tilde f_x(\tilde y),d(\Theta_{f,i}(\gamma))_{\tilde f_x(\tilde y)}\cdot\tilde v])$.
Let $V_{i,r} = C^r(\tilde C_i,\R^m)$ for $i = 1,\ldots,
N$ and $0\le r\le\infty$ with topology of uniform convergence of derivatives of order $\le r$.  This is a Banach space for finite $r$
and a \Frechet space for $r=\infty$. For finite $r$, let $\|\cdot\|_{i,r}$ be a $C^r$ norm on $V_{i,r}$.  Define a linear map
$L:\mathscr{D}^r_{\Orb}(\starfunc f^*(T\orbify{P}))\to\bigoplus_{i = 1}^N V_{i,r}$ by
\begin{equation*}
L(\sigma) = \left(\text{pr}_2(\tilde\chi_1\tilde\sigma), \ldots, \text{pr}_2(\tilde\chi_N\tilde\sigma)\right)
\end{equation*}
where $\chi_i\in C_{\Orb}^r(\orbify{O},[0,1])$, $i = 1,\ldots,N$, is
a partition of unity subordinate to the cover $\mathcal{C}$ (\cite{BB2006}*{proposition~6.1}) and $\text{pr}_2:\tilde C_i\times\R^m\to\R^m$ is bundle projection onto the second factor. Continuity of $L$ is immediate from the definitions of the $C^r$
topology on $\mathscr{D}^r_{\Orb}(\starfunc f^*(T\orbify{P}))$ and the topology on
$\bigoplus_{i =1}^N V_{i,r}$.  Moreover, given a neighborhood of the
zero section $\mathbf{0}\in\mathscr{D}^r_{\Orb}(\starfunc f^*(T\orbify{P}))$ of the
form $\mathscr{N}^r(\mathbf{0},\varepsilon_i;\mathcal{C})$, it is
apparent that there is a neighborhood of the zero section
$\mathbf{0}$ in $\bigoplus_{i=1}^N V_{i,r}$ of the form
$\max\{\|s_1\|_{1,r},\ldots,\|s_N\|_{N,r}\} < \delta$ where $\delta
\le \min\{ \varepsilon_1,\ldots,\varepsilon_N\}$ contained in
$L\left(\mathscr{N}^r(\mathbf{0},\varepsilon_i;\mathcal{C})\right)$.
Thus, with the subspace topology on
$L(\mathscr{D}^r_{\Orb}(\starfunc f^*(T\orbify{P})))$,
$L:\mathscr{D}^r_{\Orb}(\starfunc f^*(T\orbify{P}))\to L(\mathscr{D}^r_{\Orb}(\starfunc f^*(T\orbify{P})))$ is a linear homeomorphism.
Since $\mathscr{D}^r_{\Orb}(\starfunc f^*(T\orbify{P}))\subset
C^r_{\Orb}(\orbify{O},\starfunc f^*(T\orbify{P}))$ is a closed subset, we see that
$L(\mathscr{D}^r_{\Orb}(\starfunc f^*(T\orbify{P})))$ is a closed subspace of the
direct sum and thus $\mathscr{D}^r_{\Orb}(\starfunc f^*(T\orbify{P}))$ inherits a
Banach space structure if $r <\infty$ and a \Frechet space structure
if $r=\infty$.
\end{sloppypar}
\end{proof}

The following is the analogue of the notion of admissible tangent vector as defined in \cite{BB2006}.

\begin{definition}\label{AdmissibleVectors} Let $\orbify{O}, \orbify{P}$ and $f$ be as in proposition~\ref{OrbisectionsOfPullBackIsBanachSpace}. Let
  $x\in\orbify{O}$. Denote by $A_{x}(\starfunc f^*(T\orbify{P}))$ the set of {\em
    admissible vectors at $x$}
$$A_{x}(\starfunc f^*(T\orbify{P}))=\left\{v\in \starfunc f^*(T\orbify{P})_x\mid (x,v)=%
  \sigma(x)\textup{ for some
  }\sigma\in\mathscr{D}^0_{\Orb}(\starfunc f^*(T\orbify{P}))\right\}$$ 
By proposition~\ref{OrbiSecsAreVecSpace},
$A_{x}(\starfunc f^*(T\orbify{P}))$ is a vector space for each $x$, and a suborbifold of
$\starfunc f^*(T\orbify{P})_x$.  The {\em admissible pullback bundle of $T\orbify{P}$}
is the subset
$A(\starfunc f^*(T\orbify{P}))=\bigcup_{x\in\orbify{O}}A_x(\starfunc f^*(T\orbify{P}))\subset \starfunc f^*(T\orbify{P})$
with the subspace topology.  In general, $A(\starfunc f^*(T\orbify{P}))$ will not be an orbifold.
Recall that the set of admissible {\em tangent} vectors at $z$, $A_z(\orbify{P})$ as defined in \cite{BB2006} are obtained from definition~\ref{AdmissibleVectors} by replacing $\starfunc f^*(T\orbify{P})$ by $T\orbify{P}$.
\end{definition}

\section{The exponential map and proof of Theorem~\ref{MainTheorem}}\label{ProofOfMainTheoremSection}

In this section we will need several facts about Riemannian orbifolds and the exponential map. Our reference for this material will be \cite{BB2006}. Throughout this section, we assume that $\orbify{O}, \orbify{P}$ are smooth
orbifolds and that $\orbify{O}$ is compact (without boundary). Without loss of generality, we
may assume, by \cite{BB2006}*{propositions~3.11, 6.4}, that both $\orbify{O}$ and $\orbify{P}$ are
$C^\infty$ orbifolds with $C^\infty$ Riemannian metrics.

\subsection*{The exponential map} Recall the construction of the exponential map for a smooth Riemannian
orbifold $\orbify{P}$ \cite{BB2006}*{section~6}.  Assume that the
collection $\{V_\alpha\}$ is a locally finite open covering of
$\orbify{P}$ by orbifold charts that are relatively compact.  Let
$TV_\alpha \cong (\tilde V_\alpha\times\R^m)/\Gamma_\alpha$ be a local
trivialization of the tangent bundle over $V_\alpha$.  Denote the
Riemannian exponential map on $T\tilde V_\alpha$ by
$\widetilde{\exp}_{\tilde V_\alpha}:T\tilde V_\alpha\to \tilde
V_\alpha$.  Thus, for $\tilde z\in\tilde V_\alpha$ and $\tilde v\in
T_{\tilde z}\tilde V_\alpha$ we have $\widetilde{\exp}_{\tilde
  V_\alpha}(\tilde z,t\tilde v)=\tilde c_{\tilde z,\tilde v}(t)$ where
$\tilde c_{\tilde z,\tilde v}$ is the unit speed geodesic in $\tilde
V_\alpha$ which starts at $\tilde z$ and has initial velocity $\tilde
v$.  Recall that there is an open neighborhood $\tilde\Omega_{\tilde
  V_\alpha}\subset T\tilde V_\alpha$ of the $0$-section of $T\tilde
V_\alpha$ such that $\tilde c_{\tilde z,\tilde v}(1)$ is defined for
$\tilde v\in T_{\tilde z}\tilde V_\alpha\cap\tilde\Omega_{\tilde
  V_\alpha}$. Furthermore, by shrinking $\tilde\Omega_{\tilde
  V_\alpha}$ if necessary, we may assume that on $T_{\tilde z}\tilde
V_\alpha\cap\tilde\Omega_{\tilde V_\alpha}$, $\widetilde{\exp}_{\tilde
  V_\alpha}(\tilde z,\cdot)$ is a local diffeomorphism onto a
neighborhood of $\tilde z\in\tilde V_\alpha$ for each $\tilde
z\in\tilde V_\alpha$.  Let
$\Omega_\alpha=\Pi_\alpha(\tilde\Omega_{\tilde V_\alpha})$, an open
subset of $T\orbify{P}$, and define
$\Omega=\bigcup_{\alpha}\Omega_\alpha$. $\Omega$ is an open
neighborhood of the $0$-orbisection of $T\orbify{P}$.

\begin{definition}\label{ExponentialMap} Let $z\in V_\alpha$, and
  $(z,v)\in\Omega_\alpha$.  Choose $(\tilde z,\tilde v)\in\Pi_\alpha^{-1}(z,v)$. Then the {\em Riemannian exponential map}
   $\exp:\Omega\subset T\orbify{P}\to\orbify{P}$ is defined by
  $\exp(z,v)=\pi_\alpha\circ\widetilde{\exp}_{\tilde V_\alpha}(\tilde z,\tilde v)$.
\end{definition}

By \cite{BB2006}*{proposition~6.7}, this exponential map is well-defined and $\widetilde{\exp}_{\tilde V_\alpha}$ satisfies, for all 
$\delta\in\Gamma_\alpha$, the equivariance relation:
\begin{equation*}
\widetilde{\exp}_{\tilde V_\alpha}[\delta\cdot(\tilde z,\tilde v)]=\delta\cdot\widetilde{\exp}_{\tilde V_\alpha}(\tilde z,\tilde v).
\end{equation*}

As usual we denote by $\exp_z$ the restriction of $\exp$ to a single
tangent cone $T_z\orbify{P}$. We let $B(x,r)$ denote the metric
$r$-ball centered at $z$ and use tildes to denote corresponding points
in local coverings. 

\subsection*{The relation between orbisections, the exponential map and complete orbifold maps}

The composition of the exponential map with an orbisection of the pullback tangent orbibundle via a complete orbifold map $\starfunc f$ turns out
to be a smooth complete orbifold map with the same equivariance relation as $\starfunc f$.

\begin{proposition}\label{ExpOrbisectionIsMap}
  Let $\orbify{O}, \orbify{P}$ be smooth Riemannian orbifolds and let $\starfunc f\in\COrbMaps^{r+1}(\orbify{O},\orbify{P})$. Let $\sigma$ be
  a $C^r$ orbisection of the pullback tangent orbibundle $\starfunc f^*(T\orbify{P})$ with $F\circ\sigma(x)\in\Omega$. Then the map
  $E^\sigma(x)=(\exp\circ\,F\circ\,\sigma)(x):\orbify{O}\to\orbify{P}$ is a
  complete $C^r$ orbifold map with representation 
  $\starfunc E^\sigma=(E^\sigma,\{\tilde E^\sigma_x\},\{\Theta_{E^\sigma,x}\})$ where 
  $\tilde E^\sigma_x=\widetilde{\exp}_{V_{f(x)}}\circ\tilde F_x\circ\tilde\sigma_x$ and $\Theta_{E^\sigma,x}=\Theta_{f,x}$ for all $x\in\orbify{O}$.
  In particular, $\starfunc E^\sigma\in\COrbMaps^{r}(\orbify{O},\orbify{P})$.
\end{proposition}

\begin{proof} Let $(\tilde U_x,\Gamma_x)$ be an orbifold chart at
  $x\in\orbify{O}$.  For $y\in U_x$, $\sigma(y)=(y,\xi(y))\in\starfunc f^*(TV_{f(x)})$. Let $\tilde\sigma_x=(\tilde y,\tilde\xi(\tilde y))$ be a lift
  of $\sigma$ and $\tilde F_x$ the map defined after the diagram of definition~\ref{PullbackBundleDef}. Then
  the map $\tilde E^\sigma_x=\widetilde{\exp}_{V_{f(x)}}\circ\tilde F_x\circ\tilde\sigma_x$ is a $C^r$ lift of $E^\sigma$ and using the equivariance relations for $\sigma_x$, $\tilde F_x$ and $\widetilde{\exp}_{V_{f(x)}}$ we have for all $\gamma\in\Gamma_x$:
  \begin{align*}
    \tilde E^\sigma_x(\gamma\cdot\tilde y) & = \widetilde{\exp}_{\tilde V_{f(x)}}\circ\tilde F_x\circ\tilde\sigma_x(\gamma\cdot\tilde y)\\
    & = \widetilde{\exp}_{\tilde V_{f(x)}}\circ\tilde F_x(\gamma\cdot\tilde\sigma_x(\tilde y))\\
    & = \widetilde{\exp}_{\tilde V_{f(x)}}(\Theta_{f,x}(\gamma)\cdot\tilde\xi(\tilde y))\\
    & = \Theta_{f,x}(\gamma)\cdot\widetilde{\exp}_{\tilde V_{f(x)}}(\tilde\xi(\tilde y))\\
    & = \Theta_{f,x}(\gamma)\cdot\widetilde{\exp}_{\tilde V_{f(x)}}\circ\tilde F_x\circ\tilde\sigma_x(\tilde y)\\
    & = \Theta_{f,x}(\gamma)\cdot \tilde E^\sigma_x(\tilde y).
  \end{align*}
\end{proof}

\subsection*{The local manifold structure}

Let $\orbify{O}, \orbify{P}$ and $\starfunc f$ be as in proposition~\ref{ExpOrbisectionIsMap}. Denote by 
$\mathbf{0}:\orbify{O}\to \starfunc f^*(T\orbify{P})$,
$\mathbf{0}(x)=0_x\in \starfunc f^*(T\orbify{P})_x$, the $0$-orbisection of
$\starfunc f^*(T\orbify{P})$. Then $E^\mathbf{0}(x)=(\exp\circ\,F\circ\,\mathbf{0})(x)=f(x)$.
We let $\mathscr{B}_{\starfunc f}^r(\sigma,\varepsilon)=\mathscr{N}^r(\sigma,\varepsilon)\cap\mathscr{D}_{\Orb}^r(\starfunc f^*(T\orbify{P}))$.
That is, $\mathscr{B}_{\starfunc f}^r(\sigma,\varepsilon)$ is the set of $C^r$
orbisections $\varepsilon$-close to $\sigma$ in the $C^r$ topology on
$C^r_{\Orb}(\orbify{O},\starfunc f^*(T\orbify{P}))$. Proposition~\ref{ExpOrbisectionIsMap} and definition~\ref{CrTopOnComplete} immediately yield the following

\begin{proposition}\label{EIsDefined} Let $\orbify{O}, \orbify{P}$ and $\starfunc f$ be as in proposition~\ref{ExpOrbisectionIsMap} with $\orbify{O}$ compact. 
There exists $\varepsilon >0$ and continuous map $E:\mathscr{B}_{\starfunc f}^r(\mathbf{0},\varepsilon)\to\mathscr{N}^r(\starfunc f,\varepsilon)$ defined by
$E(\sigma)=E^\sigma$.
\end{proposition}

Since $\mathscr{B}_{\starfunc f}^r(\mathbf{0},\varepsilon)$ is an open subset of a Banach/\Frechet space by proposition~\ref{OrbisectionsOfPullBackIsBanachSpace}, the proof theorem~\ref{MainTheorem} will be complete if $E$ is shown to be a homeomorphism. We first show that $E$ is injective.

\begin{proposition}\label{EIsInjective} The map
  $E:\mathscr{B}_{\starfunc f}^r(\mathbf{0},\varepsilon)\to\mathscr{N}^r(\starfunc f,\varepsilon)$
  is injective.
\end{proposition}
\begin{proof} Suppose $E(\sigma)=E(\tau)$ for
  $\sigma,\tau\in\mathscr{B}_{\starfunc f}^r(\mathbf{0},\varepsilon)$. Since these are to be considered equal as complete orbifold maps, in
  each orbifold chart $(\tilde U_x,\Gamma_x)$, we must have equal local lifts:
  $\widetilde{\exp}_{\tilde V_{f(x)}}\circ\tilde F_x\circ\tilde\sigma_x(\tilde y)=%
  \widetilde{\exp}_{\tilde V_{f(x)}}\circ\tilde F_x\circ\tilde\tau_x(\tilde y)$. If we write in local coordinates
  $\tilde\sigma_x(\tilde y)=(\tilde y,\tilde\xi(\tilde y))$ and 
  $\tilde\tau_x(\tilde y)=(\tilde y,\tilde\eta(\tilde y))$ where $\xi(\tilde y)=[\tilde f_x(\tilde y),v_{\tilde f_x(\tilde y)}]$ and 
  $\eta(\tilde y)=[\tilde f_x(\tilde y),w_{\tilde f_x(\tilde y)}]$ then
 \begin{align*}
  \widetilde{\exp}_{\tilde V_{f(x)}}\circ\tilde F_x\circ\tilde\sigma_x(\tilde y)&=%
  \widetilde{\exp}_{\tilde V_{f(x)}}\circ\tilde F_x\circ\tilde\tau_x(\tilde y)\Longleftrightarrow\\
  \widetilde{\exp}_{\tilde V_{f(x)}}\circ\tilde F_x(\tilde y,\tilde\xi(\tilde y))&=%
  \widetilde{\exp}_{\tilde V_{f(x)}}\circ\tilde F_x(\tilde y,\tilde\eta(\tilde y))\Longleftrightarrow\\
  \widetilde{\exp}_{\tilde V_{f(x)}}[\tilde f_x(\tilde y),v_{\tilde f_x(\tilde y)}]&=%
  \widetilde{\exp}_{\tilde V_{f(x)}}[\tilde f_x(\tilde y),w_{\tilde f_x(\tilde y)}]
\end{align*}
  Since $\widetilde{\exp}_{\tilde V_{f(x)}}(\tilde f_x(\tilde y),\cdot)$ is a local
  $C^r$ diffeomorphism we must have $v_{\tilde f_x(\tilde y)}=w_{\tilde f_x(\tilde y)}$. Hence $\sigma=\tau$
  (as orbifold maps) and $E$ is injective.
\end{proof}

The proof of the following proposition is a slightly modified version of \cite{BB2006}*{proposition~7.3}.

\begin{proposition}\label{EIsSurjective} The map 
  $E:\mathscr{B}_{\starfunc f}^r(\mathbf{0},\varepsilon)\to\mathscr{N}^r(\starfunc f,\varepsilon)$
  is surjective.
\end{proposition}
\begin{proof}
  Let $\starfunc g=(g,\{\tilde g_x\},\{\Theta_{g,x}\})\in\mathscr{N}^r(\starfunc f,\varepsilon)$.
  Let $\{C_i\}$ be a finite covering of $\orbify{O}$ by compact sets
  such that $C_i$ is an orbifold chart and $g(C_i)\subset V_i$ where
  $V_i$ is a relatively compact orbifold chart of $\orbify{P}$. Let $x\in C_i$, and
  $\tilde U_x\subset\intr{\tilde C_i}$ an orbifold chart at $x$ where
  the local lift $\tilde g_x$ to $\tilde U_x$ is $C^0$
  $\varepsilon$-close to the local lift $\tilde f_x$. By lemma~\ref{StrataImplicationsOfLocalNeighborhood} and its proof we have
  $\Theta_{f,x}=\theta_{g(x)f(x)}\circ\Theta_{g,x}:\Gamma_x\to\Gamma_{f(x)}$. In particular, the action of $\Theta_{f,x}$ is the same as action
  $\Theta_{g,x}$ on the image $\tilde g_x(\tilde U_x)\subset\tilde V_{f(x)}$.

We wish to define a $C^r$ orbisection $\sigma$ so that $E(\sigma)=\starfunc g$.
We do this by defining appropriate local lifts $\tilde\sigma_x$. In
particular, let
$$\tilde\sigma_x(\tilde y)=(\tilde y,\tilde\xi(\tilde y))=%
\left(\tilde y,\left[\tilde f_x(\tilde y),\widetilde{\exp}^{-1}_{\tilde V_{f(x)},\tilde f_x(\tilde y)}\left(\tilde
    g_x(\tilde y)\right)\right]\right)\in \starfunc f^*(T\tilde V_{f(x)}).$$
With this definition, we see that
\begin{align*}
\tilde E^\sigma_x(\tilde y)&=\widetilde{\exp}_{\tilde V_{f(x)}}\circ\tilde F_x\circ\tilde\sigma_x(\tilde y)\\
&=\widetilde{\exp}_{\tilde V_{f(x)}}\circ\widetilde{\exp}^{-1}_{\tilde V_{f(x)},\tilde f_x(\tilde y)}(\tilde g_x(\tilde y))=\tilde g_x(\tilde y).
\end{align*}
This shows that $E(\sigma)=\starfunc g$.
All that remains to show is that $\tilde\sigma_x$ satisfies the correct equivariance relation for an orbisection. Before we do that, observe
that, in general, for $\delta\in\Gamma_{f(x)}$ we have (essentially for any exponential map)
\begin{align*}
\widetilde{\exp}^{-1}_{\tilde V_{f(x)},\tilde f_x(\tilde y)}\left(\delta\cdot\tilde z\right)&=
(d\delta)_{\delta^{-1}\tilde f_x(\tilde y)}\circ\widetilde{\exp}^{-1}_{\tilde V_{f(x)},\delta^{-1}\tilde f_x(\tilde y)}(\tilde z)\\
&=\delta\cdot\widetilde{\exp}^{-1}_{\tilde V_{f(x)},\delta^{-1}\tilde f_x(\tilde y)}(\tilde z).
\end{align*}
Thus,
\begin{align*}
  \tilde\sigma_x(\gamma\cdot\tilde y) & = (\gamma\cdot\tilde y,\tilde\xi(\gamma\cdot\tilde y)) =
  \left(\gamma\cdot\tilde y,\left[\tilde f_x(\gamma\cdot\tilde y),\widetilde{\exp}^{-1}_{\tilde V_{f(x)},\tilde f_x(\gamma\cdot\tilde y)}\left(\tilde
 g_x(\gamma\cdot\tilde y)\right)\right]\right)\\
 & = \left(\gamma\cdot\tilde y,\left[\Theta_{f,x}(\gamma)\cdot\tilde f_x(\tilde y),%
 \widetilde{\exp}^{-1}_{\tilde V_{f(x)},\Theta_{f,x}(\gamma)\tilde f_x(\tilde y)}\left(\Theta_{g,x}(\gamma)\cdot\tilde g_x(\tilde y)\right)\right]\right)\\
 & = \left(\gamma\cdot\tilde y,\left[\Theta_{f,x}(\gamma)\cdot\tilde f_x(\tilde y),%
 \Theta_{f,x}(\gamma)\cdot\widetilde{\exp}^{-1}_{\tilde V_{f(x)},(\Theta_{f,x}(\gamma))^{-1}\Theta_{f,x}(\gamma)\tilde f_x(\tilde y)}%
 \left(\tilde g_x(\tilde y)\right)\right]\right)\\
  & = \left(\gamma\cdot\tilde y,\left[\Theta_{f,x}(\gamma)\cdot\tilde f_x(\tilde y),%
 \Theta_{f,x}(\gamma)\cdot\widetilde{\exp}^{-1}_{\tilde V_{f(x)},\tilde f_x(\tilde y)}\left(\tilde g_x(\tilde y)\right)\right]\right)\\
 & = \left(\gamma\cdot\tilde y,\Theta_{f,x}(\gamma)\cdot\tilde\xi(\tilde y)\right)\\
 & = \gamma\cdot\tilde\sigma_x(\tilde y)
\end{align*}
which is the correct equivariance relation for an orbisection. As a
result we see that the map $\sigma(x)=\Pi_x\circ\tilde\sigma_x(\tilde x)$ defines a $C^r$ orbisection of $\starfunc f^*(T\orbify{P})$ and that
$E(\sigma)=\starfunc g$.
\end{proof}

The following proposition is the last step to complete the proof
of the theorem~\ref{MainTheorem}. It gives a $C^0$ manifold structure to
$\COrbMaps^r(\orbify{O},\orbify{P})$ where the model space for a neighborhood of $\starfunc f$ is the topological
vector space of $C^r$ orbisections of a pullback tangent orbibundle of $\orbify{P}$ via $\starfunc f$ with the
$C^r$ topology.

\begin{proposition}\label{EIsAHomeomorphism} The map
  $E:\mathscr{B}_{\starfunc f}^r(\mathbf{0},\varepsilon)\to\mathscr{N}^r(\starfunc f,\varepsilon)$
  is a homeomorphism.
\end{proposition}
\begin{proof} Propositions ~\ref{EIsInjective} and \ref{EIsSurjective}
  show that $E$ is bijective. Continuity of $E$ follows from the
  formula for a local lift of $E$ given in
  propositon~\ref{ExpOrbisectionIsMap} and continuity of $E^{-1}$
  follows from the formula for $\tilde\sigma_x$ given in the proof
  of proposition~\ref{EIsSurjective}.
\end{proof}

\section{Building Stratified Neighborhoods}\label{StratifiedNeighborhoodsSection}

Our first task of this section will be to prove corollary~\ref{OrbifoldMapStructure}. Let $f\in\OrbMaps^r(\orbify{O},\orbify{P})$. From the observation following definition~\ref{CrTopOnComplete}, we have that 
$$q^{-1}(\mathscr{N}^r(f,\varepsilon))=\mathscr{N}^r(\starfunc f_1,\varepsilon)\cup\cdots\cup\mathscr{N}^r(\starfunc f_k,\varepsilon)$$
is a disjoint union of neighborhoods where each complete map $\starfunc f_i=(f,\{\tilde f_x\},\{\Theta_{f,x}\}_i)$. 
We first partition the neighborhood $\mathscr{N}^r(f,\varepsilon)$. For each 
$g\in\mathscr{N}^r(f,\varepsilon)$, define $J_g\subset\{1,\ldots,k\}$ to be the set of indices $j$ such that
$q^{-1}(g)\cap\mathscr{N}^r(\starfunc f_j,\varepsilon)\ne\emptyset$. For $J\subset\{1,\ldots,k\}$, define
$$\orbify{X}_J=\{g\in\mathscr{N}^r(f,\varepsilon)\mid J=J_g\}.$$
This is a partition of $\mathscr{N}^r(f,\varepsilon)$. Of course, the partial ordering is from set inclusion: $J'\prec J\Leftrightarrow J'\subset J$.
We now verify conditions (\ref{StrataCondition1}) and (\ref{StrataCondition2}) of 
definition~\ref{StratificationDef} in the next two lemmas.

\begin{lemma}\label{XJIsAManifold} Each $\orbify{X}_J$ is a submanifold of $\mathscr{N}^r(f,\varepsilon)$.
\end{lemma}
\begin{proof} Let $J=\{j_1,\ldots,j_l\}$. For any $g\in\orbify{X}_J$, we have
$$q^{-1}(g)\cap\mathscr{N}^r(\starfunc f_{j_1},\varepsilon)\cap\cdots\cap\mathscr{N}^r(\starfunc f_{j_l},\varepsilon)\ne\emptyset.$$
By proposition~\ref{EIsAHomeomorphism}, there exists unique orbisections $\sigma_{j_1},\ldots,\sigma_{j_l}$ (of the respective pullback tangent orbibundles $\starfunc f_{j_i}^*(T\orbify{P}))$ such that $q(E(\sigma_{j_i}))=g$ for $i=1,\ldots, l$. Let $(\tilde U_x,\Gamma_x)$ be a local chart about $x\in\orbify{O}$ and let
$\tilde\sigma_{j_i,x}$ denote the local lifts of the orbisection $\sigma_{j_i}$ and let 
$\tilde F_{j_i,x}:\tilde U_x\times_{\tilde V_{f(x)}}T\tilde V_{f(x)}\to T\tilde V_{f(x)}$ denote the linear isomorphism
$(\tilde y,\tilde\xi)\mapsto\tilde\xi$ given in definition~\ref{PullbackBundleDef}
of the pullback tangent orbibundle. Since 
$\widetilde{\exp}_{\tilde V_{f(x)}}\circ\tilde F_{j_i,x}\circ\tilde\sigma_{j_i,x}(\tilde y)=\tilde g_x(\tilde y)$ for all $i=1,\ldots, l$ and since
$\widetilde{\exp}_{\tilde V_{f(x)}}(\tilde f_x(\tilde y),\cdot)$ is a local diffeomorphism, we must have
\begin{equation*}
\tilde F_{j_1,x}\circ\tilde\sigma_{j_1,x}(\tilde y)=\cdots=\tilde F_{j_l,x}\circ\tilde\sigma_{j_l,x}(\tilde y)
\end{equation*}
for all $\tilde y\in\tilde U_x$. Because $\tilde F_{j_i,x}$ is a linear isomorphism, this relation is preserved under addition and scalar multiplication of local lifts of orbisections $\sigma_{j_i,x}$. From the proof of proposition~\ref{OrbiSecsAreVecSpace}, this relation descends to 
\begin{equation}\label{XJSubmanifoldCondition}\tag{*}
F_{j_1}\circ\sigma_{j_1}(y)=\cdots=F_{j_l}\circ\sigma_{j_l}(y)
\end{equation}
for $y\in U_x$.
Since $F_{j_i}$ is a linear isomorphism when restricted to the vector space of admissible vectors $A_x(\starfunc f_{j_i}^*(T\orbify{P}))$, the set of orbisections satisfying these relations is a linear submanifold of each $\mathscr{B}_{\starfunc f_{j_i}}(\mathbf{0},\varepsilon)$. From this we may conclude that each $g\in\orbify{X}_J$ has a neighborhood modeled on a linear submanifold of $\mathscr{B}_{\starfunc f_{j_i}}(\mathbf{0},\varepsilon)$, which is enough to prove that $\orbify{X}_J$ is a submanifold of $\mathscr{N}^r(f,\varepsilon)$.
\end{proof}

\begin{lemma}\label{XJCondition2} If $\orbify{X}_J\cap\overline{\orbify{X}_{J'}}\ne\emptyset$, $J\ne J'$, then $J'\prec J$ and $\orbify{X}_J\subset\overline{\orbify{X}_{J'}}$.
\end{lemma}
\begin{proof} Let $J'=\{j_1,\ldots,j_l\}$. For $i=1,\ldots, l$, suppose $\left\{\sigma_{j_i}^{(k)}\right\}_{k=1}^\infty$ is a sequence of orbisections which converges to 
$\overline{\sigma_{j_i}}\in\mathscr{D}_{\Orb}^r(\starfunc f_{j_i}^*(T\orbify{P}))$. Further suppose each 
$\sigma_{j_i}^{(k)}$ satisfies condition \eqref{XJSubmanifoldCondition} of lemma~\ref{XJIsAManifold}. Then, by continuity, each 
$\overline{\sigma_{j_i}}$ satisfies \eqref{XJSubmanifoldCondition} also.  If we let $q(E(\overline{\sigma_{j_i}}))=\overline{g}$, and
$J=J_{\overline{g}}$ we have shown that if $\orbify{X}_J\cap\overline{\orbify{X}_{J'}}\ne\emptyset$, $J\ne J'$, then $J'\prec J$ and
$\orbify{X}_J\subset\overline{\orbify{X}_{J'}}$.
\end{proof}

Theorem~\ref{MainTheorem} with lemmas \ref{XJIsAManifold} and \ref{XJCondition2} together prove corollary~\ref{OrbifoldMapStructure}. Finally the proof of corollary~\ref{ReducedOrbifoldMapStructure} follows from corollary~\ref{OrbifoldMapStructure} and lemmas~\ref{PointwiseOrbQuotientRedOrb} and \ref{IdCompatCrTopology}. That is,
$\mathscr{N}^r(\redfunc f,\varepsilon)$ is the quotient of the finite group $\mathscr{ID}_{\orbify{N}}$ acting on 
$\mathscr{N}^r(f,\varepsilon)$. That the corresponding quotient map $q_\bullet$ restricts on each stratum to give a smooth orbifold chart follows from an argument almost identical to the argument in the proof of corollary~\ref{CompleteReducedOrbifoldMapStructure} from section~\ref{FunctionSpaceTopologySection} that $q_{\ssslozenge}$ defined a smooth orbifold chart.

\subsection*{An alternative view of the stratification}\label{GlobalView}
Up to this point, the notion of pullback bundle (definition~\ref{PullbackBundleDef}) required the use of a \emph{complete} orbifold map. 
Although not necessary for our results, we present a more global view of the stratification obtained above by defining directly an appropriate notion of pullback for an orbifold map $f$.  We will use the setup of this section and the notation of definition~\ref{PullbackBundleDef}. However, for convenience we will write a complete orbifold map 
$\starfunc f_i=(f,\{\tilde f_x\},\{\Theta_{f,x,i}\})$.
To begin, we let $f^*(T\orbify{P})$ be the
space defined by:
\begin{equation*}
 f^*(T\orbify{P}) = \left.\left(\bigsqcup_{i = 1}^k \starfunc f_i^*(T\orbify{P})\right)\right/\sim
\end{equation*}
where the equivalence relation $\sim$ is defined as follows: Let 
$$\bigsqcup_{i=1}^k\left(\tilde U_x\times_{\tilde V_{f(x)}}T\tilde V_{f(x)}\right)_i$$
denote the disjoint union of $k$ copies of $\tilde U_x\times_{\tilde V_{f(x)}}T\tilde V_{f(x)}$. Then in local bundle charts, for
$(\tilde y_i, \tilde \xi_i)\in\left(\tilde U_x\times_{\tilde V_{f(x)}}T\tilde V_{f(x)}\right)_i$,

\begin{gather*}
(\tilde y_i, \tilde \xi_i) \sim  (\tilde y_j, \tilde \xi_j) \Longleftrightarrow \\
\tilde y_i = \tilde y_j\text{ and,}\\
\Theta_{f, x, i}(\gamma)\cdot\tilde \xi_i =\Theta_{f, x, j}(\gamma) \cdot\tilde \xi_j\text{ for all }\gamma \in \Gamma_x.
\end{gather*}
There is an obvious projection map onto $\orbify{O}$ and the total space of
$f^*(T\orbify{P})$ is a bundle over $\orbify{O}$. Note that there are standard continuous injections
$\iota_i:\starfunc f_i^*(T\orbify{P})\to f^*(T\orbify{P})$ and that the bundle maps $F_i:\starfunc
f_i^*(T\orbify{P})\to T\orbify{P}$ glue together to give a continuous
bundle map $F:f^*(T\orbify{P})\to T\orbify{P}$ satisfying $F\circ\iota_i = F_i$.

We also define for $J\subset \{1, \ldots, k\}$ the suborbifold $f^*(T\orbify{P})_J$ of $f^*(T\orbify{P})$  by
\begin{equation*}
 f^*(T\orbify{P})_J = \bigsqcup_{i \in J} \starfunc f_i^*(T\orbify{P}_J)/\sim
\end{equation*}
where $T\orbify{P}_J$ is the subspace of $T\orbify{P}$
covered in bundle charts by
$$(T\tilde V_{f(x)})_J=\left\{(\tilde f_x(\tilde y),\tilde \xi)\in T\tilde V_{f(x)}\mid
\Theta_{f, x, i}(\gamma)\cdot\tilde \xi = \Theta_{f, x, j}(\gamma)\cdot\tilde\xi\text{ for all }i, j \in J\right\}.$$
Finally, let
\begin{multline*}
 \mathscr{D}^r_{\Orb}(\starfunc f_i^*(T\orbify{P}))_J = \\
 \{\sigma\in\mathscr{D}^r_{\Orb}(\starfunc f_i^*(T\orbify{P}))\mid \tilde F_{x,i}\circ \tilde\sigma_x(\tilde y) \in (T\tilde V_{f(x)})_J \text{ for all } \tilde y \in \tilde U_x \text{ and }x \in \orbify{O} \}.
\end{multline*}
Note that $\tilde F_{x,i}\circ\tilde F_{x,j}^{-1}:\mathscr{D}^r_{\Orb}(\starfunc f_j^*(T\orbify{P}))_J\to
\mathscr{D}^r_{\Orb}(\starfunc f_i^*(T\orbify{P}))_J$ is a linear
isomorphism of Banach ($r$ finite)/\Frechet ($r=\infty$) subspaces for all $i$,
$j\in J$ and $0\le r\le \infty$.  By abuse of notation, we write
$\mathscr{D}^r_{\Orb}(f^*(T\orbify{P}))$ for the space of $C^r$ orbisections
$\sigma:\orbify{O}\to f^*(T\orbify{P})$ equipped with the $C^r$
topology.  From the construction of $f^*(T\orbify{P})$ it is clear that the
Riemannian exponential map on $T\orbify{P}$ induces a map
$E:\mathscr{D}_{\Orb}^r(f^*T\orbify{P})\to \OrbMaps^r(\orbify{O},\orbify{P})$ as in proposition~\ref{EIsDefined}.  For $f\in \OrbMaps^r(\orbify{O}, \orbify{P})$ and
$\starfunc f_i\in \COrbMaps^r(\orbify{O}, \orbify{P})$ mapping to $f$
we let $\Theta(f)_x = \{\Theta_{f,x,i}\}$ where $\starfunc{f}_i = (f,\{\tilde f_x\},\{\Theta_{f,x,i}\})$.
\begin{lemma}
 There are neighborhoods $\mathscr{B}_f^r(\mathbf{0},\varepsilon)$ of $\mathbf{0}\in \mathscr{D}_{\Orb}^r(f^*T\orbify{P})$ and 
 $\mathscr{N}^r(f,\varepsilon)$ of $f$ in $\OrbMaps^r(\orbify{O},\orbify{P})$  so that
 $E:\mathscr{B}_f^r(\mathbf{0},\varepsilon)\to \mathscr{N}^r(f,\varepsilon)$ is a homeomorphism.
\end{lemma}
\begin{proof}
 The proof follows from observing that $\{g\in \OrbMaps^r(\orbify{O},\orbify{P}) \mid \Theta(g)_x \subset \Theta(f)_x $ for all $x\in
 \orbify{O}\}$ is an open subset (since the homomorphisms $\Theta_{f,x}$ are locally constant).  By theorem~\ref{MainTheorem}, there is a neighborhood
 of each $\starfunc{f}_i$ for which the  map $E$ of proposition~\ref{EIsDefined} is a
 homeomorphism. By taking $g\in \OrbMaps^r(\orbify{O}, \orbify{P})$
 as above and sufficiently $C^r$ close to $f$, all of of its preimages
 $\starfunc g_j = (g, \{\tilde g_x\}, \{\Theta_{g, x, j}\})$ will lie in such
 neighborhoods.
\end{proof}

\section{Some Infinite-dimensional Analysis}\label{InfiniteDimAnalysis} In this section we recall the results of global analysis that we need in order to substantiate our various claims of smoothness. For finite order differentiability, a strong argument can be made that the Lipschitz categories $\text{Lip}^r$ are better suited to questions of calculus than the more common $C^r$ category \cite{MR961256}. For our purposes, however, we have chosen to use the $C^r$ category for questions of finite order differentiability, and for questions of infinite order differentiability, we use the convenient calculus as detailed in the monographs \cites{MR961256, MR1471480}. 

\subsection*{Review of the convenient calculus}\label{ConvenientCalculus} For any topological vector space $E$ the notion of a smooth curve $c\in C^\infty(\mathbb{R},E)$ makes sense using the usual difference quotient and iterating. A mapping
$f:E\to F$ between locally convex vector spaces is called $\emph{smooth}$ if it maps smooth curves to smooth curves. That is, if $f\circ c\in C^\infty(\R,F)$ for all $c\in C^\infty(\R,E)$. For $E$, $F$ finite dimensional this yields the usual notion of ($C^\infty$-) smoothness. Unfortunately, such a characterization fails for finite order 
($C^r$-) differentiability \cite{MR0237728}. Generalizing the fact that a map $f$ between finite dimensional vector spaces is smooth if and only if its component functions are smooth, Fr\"{o}licher, Kriegl and Michor \cites{MR961256, MR1471480} introduce the notion of a \emph{convenient} vector space: A locally convex vector space is \emph{convenient} if every scalarwise smooth curve $c:\R\to E$ is smooth. $c$ is a scalarwise smooth curve if $\ell\circ c:\R\to\R$ is smooth for all continuous linear functionals $\ell$ on $E$. For our purposes, we remark that if $E$ is a \Frechet space then $E$ is convenient and the locally convex topology agrees with the $c^\infty$-topology or final topology  with respect to the set of mappings $C^\infty(\R,E)$. It then follows from these definitions that smooth mappings between \Frechet (or, convenient vector) spaces are continuous. At this point, one can work on open subsets of convenient vector spaces and introduce the notions of smooth manifold, smooth tangent bundle and smooth Lie group modeled on convenient vector spaces in a straightforward way.

An important feature of the convenient setting is the following theorem on Cartesian closedness:

\begin{theorem}[\cite{MR1471480}*{Theorem~3.12}]\label{CartesianClosedness} Let $A_i\subset E_i$ be $c^\infty$-open subsets in locally convex spaces, which need not be $c^\infty$-complete. Then a mapping $f:A_1\times A_2\to F$ is smooth if and only if the canonically associated mapping
$f^\vee:A_1\to C^\infty(A_2,F)$ exists and is smooth.
\end{theorem}

Our first use of the convenient calculus will be to substantiate the smoothness claims of Theorem~\ref{MainTheorem} completing the proof. Throughout the remainder, we assume, as in section~\ref{ProofOfMainTheoremSection}, that all orbifolds are $C^\infty$ with $C^\infty$ Riemannian metric. Further, the orbifold $\orbify{O}$ will be compact (without boundary).

\begin{lemma}\label{COrbMapsAreSmoothManifold} $\COrbMaps^r(\orbify{O},\orbify{P})$ has the structure of a smooth ($C^\infty$) Banach/\Frechet manifold.
\end{lemma}
\begin{proof} We have already shown in section~\ref{ProofOfMainTheoremSection} that $\COrbMaps^r(\orbify{O},\orbify{P})$ has the required structure as a topological Banach/\Frechet manifold. Let 
$\starfunc f\in\COrbMaps^r(\orbify{O},\orbify{P})$ and let 
$E_{\starfunc f}:\mathscr{B}_{\starfunc f}^r(\mathbf{0},\varepsilon')\to\mathscr{N}^r(\starfunc f,\varepsilon')$ be a manifold chart about $\starfunc f$ given by proposition~\ref{EIsAHomeomorphism}. Let 
$\starfunc g\in\mathscr{N}^r(\starfunc f,\varepsilon')$ and choose $0< \varepsilon$ $ < \varepsilon '$
so that the manifold chart
$E_{\starfunc g}:\mathscr{B}_{\starfunc g}^r(\mathbf{0},\varepsilon)\to\mathscr{N}^r(\starfunc g,\varepsilon)$ is contained entirely within
$\mathscr{N}^r(\starfunc f,\varepsilon')$. Then the chart transition map
$$E_{\starfunc f}^{-1}\circ E_{\starfunc g}:\mathscr{B}_{\starfunc g}^r(\mathbf{0},\varepsilon)\to%
\left(E_{\starfunc f}^{-1}\circ E_{\starfunc g}\right)\left(\mathscr{B}_{\starfunc g}^r(\mathbf{0},\varepsilon)\right)\subset%
\mathscr{B}_{\starfunc f}^r(\mathbf{0},\varepsilon')$$
is a homeomorphism between open subsets of Banach/\Frechet spaces. Thus to show smoothness, we need to show that 
$E_{\starfunc f}^{-1}\circ E_{\starfunc g}$ takes smooth curves to smooth curves. So, let 
$\sigma^t:\R\to\mathscr{B}_{\starfunc g}^r(\mathbf{0},\varepsilon)\subset\COrbMaps^r(\orbify{O},\starfunc g^*(T\orbify{P}))$ be a smooth curve and let 
$\tilde\sigma^t_x:(0,1)\to\mathscr{D}^r(\starfunc g^*(T\tilde V_{g(x)}))\subset%
C^r(\tilde U_x,\tilde U_x\times_{{\tilde V}_{g(x)}}T\tilde V_{g(x)})$ be a smooth local equivariant lift over an orbifold chart $U_x$. The interval $(0,1)$ is being chosen for convenience to make clear that we want the image of the lift $\tilde\sigma^t_x$ to lie in a single trivializing bundle chart. The key observation is that the computations of difference quotients for $\sigma^t$ are identical in the local lift $\tilde\sigma^t_x$ since an orbisection must take values in the admissible bundle $A(\starfunc g^*(T\orbify{P}))$ whose fibers are the vector spaces fixed by the action of the local isotropy subgroups $\Gamma_x$ (section~\ref{TangentOrbibundleSection}). In particular, it follows that $\sigma^t$ is smooth if and only if each local lift $\tilde\sigma^t_x$ is smooth.
Using lemma~\ref{CartesianClosedness}, it follows that
$\tilde\sigma^\wedge_x(t,\tilde y):(0,1)\times\tilde U_x\to\tilde U_x\times_{{\tilde V}_{g(x)}}T\tilde V_{g(x)}$ is smooth. From the formulas in section~\ref{ProofOfMainTheoremSection} for $E$ and its local lifts and using the 
($C^\infty$-) smoothness of the Riemannian exponential map, it follows that the map
$\tilde\eta^\wedge_x(t,\tilde y)=\tilde E_{\starfunc f,x}^{-1}\circ%
\tilde E_{\starfunc g,x}\circ\tilde\sigma^\wedge_x(t,\tilde y):(0,1)\times\tilde U_x\to%
\tilde U_x\times_{{\tilde V}_{f(x)}}T\tilde V_{f(x)}$ is smooth. Another application of lemma~\ref{CartesianClosedness} implies that $\tilde\eta^t_x=\tilde E_{\starfunc f,x}^{-1}\circ%
\tilde E_{\starfunc g,x}\circ\tilde\sigma^t_x:(0,1)\to%
\mathscr{D}^r(\starfunc f^*(T\tilde V_{f(x)}))\subset%
C^r(\tilde U_x,\tilde U_x\times_{{\tilde V}_{f(x)}}T\tilde V_{f(x)})$ is smooth. Thus, by our earlier observation,
$\eta^t=\left(E_{\starfunc f}^{-1}\circ E_{\starfunc g}\right)(\sigma^t):%
\R\to\mathscr{B}_{\starfunc f}^r(\mathbf{0},\varepsilon')\subset\COrbMaps^r(\orbify{O},\starfunc f^*(T\orbify{P}))$ is a smooth curve.
\end{proof}

A useful consequence of the observation made in lemma~\ref{COrbMapsAreSmoothManifold} is the following (compare \cite{MR1471480}*{Lemma~42.5}):
\begin{corollary}\label{CurveSmoothEquivLiftSmooth}The following conditions on a curve 
$c=(f^t,\{\tilde f^t_x\},\{\Theta_{f^t,x}\}):\R\to\COrbMaps^r(\orbify{O},\orbify{P})$ are equivalent:
\begin{enumerate}
\item\label{cSmooth1} $c$ is smooth
\item\label{cSmooth2} each local equivariant lift $\tilde f^t_x:(0,1)\to C^r(\tilde U_x,\tilde V_z)$ is smooth
\item\label{cSmooth3} each local equivariant lift $\tilde f^\wedge(t,x):(0,1)\times\tilde U_x\to\tilde V_z$ is smooth
\end{enumerate}
\end{corollary}
\begin{proof} Note that the interval $(0,1)$ is being chosen for convenience to make clear that we want the image of the lift $\tilde f^t_x$ to lie in a single orbifold chart. (\ref{cSmooth1}) $\Longleftrightarrow$ (\ref{cSmooth2}) follows from the observation that smoothness of $c$ is equivalent in local charts to smoothness of a curve into a space of orbisections. This in turn is equivalent to smoothness of the local equivariant lifts of the orbisections as in lemma~\ref{COrbMapsAreSmoothManifold}. This, of course, is equivalent to smoothness of the local equivariant lifts $\tilde f^t_x$. (\ref{cSmooth2}) $\Longleftrightarrow$ (\ref{cSmooth3}) follows from theorem~\ref{CartesianClosedness}.
\end{proof}

\begin{lemma}\label{ActionOfIdIsSmooth} Let 
$I:\mathscr{N}^s(\starfunc f,\varepsilon)\to\mathscr{N}^s(I\circ \starfunc f,\varepsilon)$, 
$\starfunc g\mapsto I\circ\starfunc g$, be the homeomorphism from lemma~\ref{IdCompatCrTopology}. Then $I$ is a $C^\infty$ diffeomorphism.
\end{lemma}
\begin{proof} By lemma~\ref{IdCompatCrTopology} and theorem~\ref{MainTheorem} via lemma~\ref{COrbMapsAreSmoothManifold}, we know that $I$ is a map between open subsets of a smooth $C^\infty$ Banach/\Frechet manifold. By
\cite{MR1471480}*{section~27}, it suffices to show that
$$\sigma^t=E^{-1}_{I\circ\starfunc f}\circ I\circ\starfunc g^t:\R\to\mathscr{D}^s_{\Orb}((I\circ\starfunc f))^*(T\orbify{P}))$$
is smooth for any smooth curve $\starfunc g^t=(g^t,\{\tilde g_x^t\},\{\Theta_{g^t,x}\})$. From an argument similar to that given in lemma~\ref{COrbMapsAreSmoothManifold}, it follows that each local equivariant lift
$\tilde g_x^t:(0,1)\to C^s(\tilde U_x,\tilde V_{f(x)})$ is smooth. Using the discussion after example~\ref{IdentityMap}, for $I=(\text{Id},\{\eta_x\cdot\tilde y\})$ we have
$$I\circ\starfunc g^t=(g^t,\{\eta_x\cdot\tilde g_x^t\},\{\gamma\mapsto\eta_x\Theta_{g^t,x}(\gamma)\eta_x^{-1}\}).$$
Then by the formulas in proposition~\ref{EIsSurjective}, $\sigma^t$ has local equivariant lifts
\begin{multline*}
\tilde\sigma^t_x(\tilde y)=(\tilde y,\tilde\xi^t(\tilde y))=\\%
\left(\tilde y,\left[\eta_x\cdot\tilde f_x(\tilde y),\widetilde{\exp}^{-1}_{\tilde V_{f(x)},\eta_x\cdot\tilde f_x(\tilde y)}\left(\eta_x\cdot\tilde g^t_x(\tilde y)\right)\right]\right)\in(I\circ \starfunc f)^*(T\tilde V_{f(x)}).
\end{multline*}
Since the exponential map and action of local isotropy subgroups are ($C^\infty$-) smooth in charts and $\tilde g_x^t$ is smooth, we see that $\tilde\sigma^t_x$ is smooth. It then follows that $\sigma^t$ is smooth by corollary~\ref{CurveSmoothEquivLiftSmooth} or the observation in the proof of lemma~\ref{COrbMapsAreSmoothManifold}, and this completes the proof.
\end{proof}

The next result we will need is that composition in our spaces of smooth orbifold maps is smooth.

\begin{lemma}\label{CompositionIsSmooth} Let $\orbify{O}$, $\orbify{P}$ and $\orbify{R}$ be smooth $C^\infty$ compact orbifolds without boundary. Then the composition mappings
\begin{align*}
\textup{comp}&: \COrbMaps^\infty(\orbify{P},\orbify{R})\times\COrbMaps^\infty(\orbify{O},\orbify{P})\to
\COrbMaps^\infty(\orbify{O},\orbify{R}),\ \ &(\starfunc f, \starfunc g)&\mapsto \starfunc f\circ \starfunc g,\\
\textup{comp}&: \CROrbMaps^\infty(\orbify{P},\orbify{R})\times\CROrbMaps^\infty(\orbify{O},\orbify{P})\to
\CROrbMaps^\infty(\orbify{O},\orbify{R}),\ \ %
&(\lozengefunc f, \lozengefunc g)&\mapsto \lozengefunc f\circ \lozengefunc g,\\
\textup{comp}&:\OrbMaps^\infty(\orbify{P},\orbify{R})\times\OrbMaps^\infty(\orbify{O},\orbify{P})\to
\OrbMaps^\infty(\orbify{O},\orbify{R}),\ \ &(f,g)&\mapsto f\circ g,\\
\textup{comp}&:\RedOrbMaps^\infty(\orbify{P},\orbify{R})\times\RedOrbMaps^\infty(\orbify{O},\orbify{P})\to
\RedOrbMaps^\infty(\orbify{O},\orbify{R}),\ \ &(\redfunc f,\redfunc g)&\mapsto \redfunc f\circ \redfunc g,\\
\end{align*}
are smooth.
\end{lemma}
\begin{proof} By lemmas~\ref{PointwiseCOrbQuotientCROrb}, \ref{PointwiseOrbQuotientRedOrb}, \ref{IdCompatCrTopology}, and the observation following definition~\ref{CrTopOnComplete}, it suffices to prove the result for the complete orbifold maps.
We use an argument analogous to \cite{MR1471480}*{Theorem~42.13}. Namely, 
let $(c_1,c_2):\R\to\COrbMaps^\infty(\orbify{P},\orbify{R})\times\COrbMaps^\infty(\orbify{O},\orbify{P})$ be a smooth curve. Then $(\text{comp}\circ(c_1,c_2))(t)(x)=c_1^\wedge(t,c_2^\wedge(t,x))$ is smooth by corollary~\ref{CurveSmoothEquivLiftSmooth}. Hence $\text{comp}$ is smooth.
\end{proof}

Finally, for finite order differentiability, we will need to refer to the $\Omega$-lemma of Palais \cite{MR0248880} as stated in \cite{MR2096566}:

\begin{lemma}[$\Omega$-lemma]\label{OmegaLemma} Let $M$ be a $C^\infty$ compact manifold and let $\tau:E\to M$ and $\tau':E'\to M$ be $C^r$ vector bundles over $M$. Let $U\subset E$ be open and let $\omega:U\subset E\to E'$ be a $C^\infty$ vector bundle map. Then the induced map
$$\Omega_\omega:\mathscr{D}^r(U)\subset\mathscr{D}^r(\tau)\to\mathscr{D}^r(\tau'),\ \Omega_\omega(\xi)=\omega\circ\xi$$
is a $C^\infty$ map. If $\omega$ is only $C^{r+k}$, then $\Omega_\omega$ is $C^{k}$.
\end{lemma}

\end{document}